\documentclass[a4paper,11pt,titlepage]{article}
\usepackage[english]{babel}%
\usepackage{epsfig}%
\usepackage{times}
\usepackage{rotating}%
\usepackage{amsmath}%
\usepackage{amssymb}
\usepackage[dvipdfm]{hyperref}
\usepackage[dvips]{color}
\def\bbeta{\mbox{\boldmath $\beta$}}
\def\bgamma{\mbox{\boldmath $\gamma$}}
\def\bdelta{\mbox{\boldmath $\delta$}} \def\bDelta{\mbox{\boldmath $\Delta$}}
\def\bxi{\mbox{\boldmath $\xi$}}
\def\bepsilon{\mbox{\boldmath $\epsilon$}}
 \def\bmu{\mbox{\boldmath $\mu$}}
 \def\bmX{\mathbf{\mathcal{X}}}
 \def\bOmega{\mbox{\boldmath $\Omega$}}
\def\bxi{\mbox{\boldmath $\xi$}}

\def\bxi{\mbox{\boldmath $\xi$}}
\def\btheta{\mbox{\boldmath $\theta$}}

\def\w{\mbox{w}} 
 
\def\bi{\mathbf{i}} \def\bA{\mathbf{A}}
\def\bE{\mathbf{E}} \def\h{\mbox{h}}
\def\bK{\mathbf{K}}
\def\bQ{\mathbf{Q}}
\def\bH{\mathbf{H}}

\def\bD{\mathbf{D}}\def\bU{\mathbf{U}}
 \def\bh{\mathbf{h}}
 
 \def\bl{\mbox{\boldmath $l$}}

\def\bL{\mathbf{L}}

\def\bw{\mathbf{w}}
\def\by{\mathbf{y}} 
\def\0{\mbox{\bf{0}}} 
\def\bS{\mathbf{S}}
\def\bI{\mathbf{I}}\def\be{\mathbf{e}}
\def\bz{\mathbf{z}}\def\bZ{\mathbf{Z}}
\def\bC{\mathbf{C}}\def\bM{\mathbf{M}}
\def\bc{\mathbf{c}} \def\bT{\mathbf{T}}\def\bm{\mathbf{m}}
\def\bW{\mathbf{W}}  
 \def\bX{\mathbf{X}} \def\bx{\mathbf{x}}
\def\bv{\mathbf{v}}

          \def\E{\mbox{E}}     \def\diag{\mbox{diag}}
      
\def\NID{\mbox{NID}} \def\N{\mbox{N}}

\RequirePackage{comment}
\excludecomment{frontmatter}

\begin{document}

\begin{frontmatter}

\title{On the Spectral Properties of Matrices Associated with
Trend Filters}

\begin{aug}
\author[A]{\fnms{Alessandra} \snm{Luati}\corref{}\ead[label=e1]{alessandra.luati@unibo.it}}
\author[B]{\fnms{Tommaso} \snm{Proietti}}
\address[A]{Dip. Scienze Statistiche\\ University of Bologna}
\address[B]{S.E.F. e ME. Q. \\
University of Rome ``Tor Vergata''}
\affiliation{University of Bologna and University of Rome}
\end{aug}

\begin{abstract}
This paper is concerned with the spectral properties of matrices
associated with linear filters for the estimation of the
underlying trend of a time series. The interest lies in the fact
that the eigenvectors can be interpreted as the latent components
of any time series that the filter smooths through the
corresponding eigenvalues. A difficulty arises because matrices
associated with trend filters are finite approximations of
Toeplitz operators and therefore very little is known about their
eigenstructure, which also depends on the boundary conditions or,
equivalently, on the filters for trend estimation at the end of
the sample.

Assuming reflecting boundary conditions, we derive a time series
decomposition in terms of periodic latent components and
corresponding smoothing eigenvalues. This decomposition depends on
the local polynomial regression estimator chosen for the interior.
Otherwise, the eigenvalue distribution is derived with an
approximation measured by the size of the perturbation that
different boundary conditions apport to the eigenvalues of
matrices belonging to algebras with known spectral properties,
such as the Circulant or the Cosine. The analytical form of the
eigenvectors is then derived with an approximation that involves
the extremes only.

A further topic investigated in the paper concerns a strategy for
a filter design in the time domain. Based on cut-off eigenvalues,
new estimators are derived, that are less variable and almost
equally biased as the original estimator, based on all the
eigenvalues. Empirical examples illustrate the effectiveness of
the method.
\end{abstract}



\end{frontmatter}

\title{\bf On the Spectral Properties of Matrices Associated with
Trend Filters}
\author{Alessandra Luati \\ Dip.
Scienze Statistiche\\ University of Bologna \and Tommaso Proietti\\ S.E.F. e ME. Q. \\
University of Rome ``Tor Vergata'' }\maketitle

\paragraph{\large{On the Spectral Properties of Matrices Associated with Trend Filters}} \vspace{2cm}

\paragraph{Abstract}

This paper is concerned with the spectral properties of matrices
associated with linear filters for the estimation of the
underlying trend of a time series. The interest lies in the fact
that the eigenvectors can be interpreted as the latent components
of any time series that the filter smooths through the
corresponding eigenvalues. A difficulty arises because matrices
associated with trend filters are finite approximations of
Toeplitz operators and therefore very little is known about their
eigenstructure, which also depends on the boundary conditions or,
equivalently, on the filters for trend estimation at the end of
the sample.

Assuming reflecting boundary conditions, we derive a time series
decomposition in terms of periodic latent components and
corresponding smoothing eigenvalues. This decomposition depends on
the local polynomial regression estimator chosen for the interior.
Otherwise, the eigenvalue distribution is derived with an
approximation measured by the size of the perturbation that
different boundary conditions apport to the eigenvalues of
matrices belonging to algebras with known spectral properties,
such as the Circulant or the Cosine. The analytical form of the
eigenvectors is then derived with an approximation that involves
the extremes only.

A further topic investigated in the paper concerns a strategy for
a filter design in the time domain. Based on cut-off eigenvalues,
new estimators are derived, that are less variable and almost
equally biased as the original estimator, based on all the
eigenvalues. Empirical examples illustrate the effectiveness of
the method.


\newpage

\section{Introduction}

The smoothing problem has a long and well established tradition in
statistics and has a wide range of applications in time series
analysis; see Anderson (1971, ch. 3), Kendall (1973), Kendall,
Stuart and Ord (1983) and Cleveland and Loader (1996). In its
simplest form, it aims at providing a measure of the underlying
tendency from noisy observations, and takes the name of signal
extraction in engineering, trend estimation in econometrics, and
graduation in actuarial sciences. This paper is concerned with
local polynomial regression methods, that developed as an
extension of least squares regression and result in estimates that
are linear combinations of the available information. These linear
combinations are often termed filters and their analysis provides
useful insight into what the method does.

The properties of linear filters are traditionally studied on
different, complementary viewpoints. In the time domain, the
analysis of the filter weights provides information on the amount
of bias introduced and variance left in the input data from the
smoothing procedure. In the frequency domain, the basic assumption
is that a time series can be decomposed as a linear combination of
trigonometric functions. The variability and the dependence
relation among the variables are then evaluated in terms of the
contribution of such components with respect to some frequency or
periodicity, usually measured in radians.

An alternative approach consists of analysing the matrices
associated with linear filters. Though smoothers have been
introduced in a time series framework, with the works of Whittaker
(1923) on spline smoothing and of Henderson (1916, 1924) on
graduation by averages, they have been mainly analysed in the
context of linear regression and in generalised additive models,
following the approach of Buja, Hastie and Tibshirani (1989,
section 2) and Hastie and Tibshirani (1990, section 3.7), based on
the smoother matrices associated with linear estimators. In these
references, the attention is concentrated on symmetric matrices
that arise as the solutions of penalised least squares problems,
such as the cubic smoothing spline estimators (see Whaba, 1990,
and Green and Silverman,
1994). 
The spectral properties of smoother matrices are analysed and
inferential procedures based on eigenvalues and eigenvectors are
developed. The authors remark that eigenanalysis is no longer
useful for non symmetric smoother matrices because of complex
eigenvalues and eigenvectors and argue that the spectral analysis
of a smoother matrix is closely related to the study of the
transfer function of the associated linear filter for time series.

These two remarks motivated the present paper. In considering
local polynomial regression methods for the estimation of the
underlying trend of a time series, symmetry is in general lost and
replaced by centrosymmetry. At the same way, the interpretation
that can be ascribed to the eigenvalues and eigenvectors of time
series smoothing matrices (let us suppose for the moment that we
are capable to lead the problem to the real or to the symmetric
case) provides useful information on the estimation method. In
fact, the eigenvectors of matrices associated with local
polynomial regression estimators can be interpreted as the latent
components of any time series that the filter smooths through the
corresponding eigenvalues. This interpretation allows a
decomposition of a time series in periodic latent components that
depend on the estimation method and opens the way to
eigenvalue-based inferential procedures. Furthermore, it is
possible to establish a formal connection between the spectrum of
a smoothing matrices and the transfer function of the associated
filter.

This paper analyses the spectral properties of matrices associated
with trend filters. In referring to spectral properties in a time
series setting, we shall distinguish between two accomplished
theories: the spectral analysis of a linear filter, where the
filter properties are studied in the frequency domain, and the
spectral properties of the associated matrix, i.e. the study of
its eigenvalues and eigenvectors. Both these techniques are
related to the concept of spectrum, to be intended as a latent
characteristic that cannot be directly observed. The spectral
properties of linear filters have been widely investigated in time
series analysis, where classical references are the books by
Jenkins and Watts (1968), Priestley (1981), Bloomfield (2000). On
the other hand, the spectral properties of the associated matrices
have not been explored. One reason is certainly due to the lack of
attention surrounding time series smoothing matrices. Another
justification relies on the fact that the mathematics of these
matrices is rather problematical. In fact, they can be interpreted
as finite approximations of infinite symmetric banded Toeplitz
operators. The latter have been extensively explored, but their
finite counterparts subject to boundary conditions are much more
difficult to analyse (see B\"{o}ttcher and Grudsky, 2005; see also
Gray, 2006). Established results hold for tridiagonal matrices,
but when the span of the filter increases, the algebra becomes
extremely complicated and, except for some cases, only approximate
results can be obtained. The size of the approximation essentially
depends on the boundary conditions on the finite matrix.
Furthermore, the boundary conditions determine the asymmetric
filters for the estimation of the trend at the extremes of the
series. Specifically, two-sided symmetric filters cannot be
applied since future (or past) observations are not available. It
should be remarked that the estimates at the end of the sample are
crucial in current analysis.

We derive approximate results on the eigenvalues and eigenvectors
of matrices associated with trend filters by interpreting the
latter as perturbations of matrices belonging to the circulant and
to the reflecting algebras, for which eigenvalues and eigenvectors
can be known exactly even in finite dimensions. The underlying
hypothesis is that of a circular and of a reflecting process,
respectively. The key result is a perturbation theorem that draws
some conclusions on the distribution of the eigenvalues of the
original smoothing matrices. We then relate the absolute
eigenvalue distribution to the gain function of the corresponding
symmetric filter. To illustrate these results, we consider a class
of asymmetric filters that approximate a given symmetric estimator
with a minimum mean square revision error strategy, subject to
polynomial constraints. This class encompasses the local
polynomial regression filters that automatically adapt at the
boundaries and that under mild assumptions on the trend are
unbiased estimators. Concerning the eigenvectors, we show that
filters that are unbiased with respect to polynomial trends of
order $p$ have $p+1$ eigenvectors that describe polynomial
functions up to the degree $p$. The analytical form of the
remaining eigenvectors is derived with an approximation which
involves the extremes only. A further topic investigated in the
paper concerns a strategy for a filter design in the time domain.
Based on cut-off eigenvalues, it is possible to obtain new
estimators that, in the interior, have less variance and almost
equal bias than the original estimator. The effectiveness of this
method is illustrated with empirical examples. We would like to
remark that even if these results are derived in a time series
setting, they apply to any non symmetric banded smoother matrix.

The paper is organised as follows. Section \ref{sec:lpr} reviews
the derivation of linear smoothers for trend extraction, both in
the interior and at the boundaries (section \ref{sec:bound}),
providing examples that will be used for the applications of the
methods developed later on in the paper. In section \ref{sec:sm},
time series smoothing matrices are introduced and their properties
are illustrated. Section \ref{sec:sa} contains the major results
of the paper, i.e. the spectrum analysis of matrices associated
with trend filters. Specifically, two sets of boundary conditions
are are considered, circulant (section \ref{sec:cbc}) and
reflecting (section \ref{sec:rbc}). Furthermore, we provide the
interpretation of the eigenvectors as analytical periodic
functions of the time. In section \ref{sec:fd}, a strategy for a
filter design based on a selected number of latent components is
derived, based on a suitably chosen cut-off eigenvalue. The
bias-variance trade off between old and new estimators is
evaluated (section \ref{sec:bv}) and the new filters are applied
to real data (section \ref{sec:ill}). Section \ref{sec:c}
summarises and comments on the results. Proofs and other technical
details are given in section \ref{sec:app}.

\section{Local polynomial regression methods}\label{sec:lpr}

Time series analysis is often based on additive models like
\begin{equation}
y_t = \mu_t +\epsilon_t, t= 1, \ldots, n,
\label{eq:y}\end{equation} where $y_t$ is the observed time
series, $\mu_t$ is the trend component, also termed the signal,
and $\epsilon_t$ is the noise, or irregular, component. The signal
$\mu_t$ can be a random or deterministic smooth function of time
whereas the most common assumption for the noise $\epsilon_t$ is
that it follows a zero mean stochastic process, such as a White
Noise or/and Gaussian. Let us assume that in (\ref{eq:y}) $\mu_t$
is an unknown deterministic function of time, so that $\E(y_t) =
\mu_t$, and that equally spaced observations $y_{t+j}, j = 0, \pm
1,2, \ldots, h,$ are available in a neighbourhood of time $t$. Our
interest lies in estimating the level of the trend at time $t$,
$\mu_t$, using the available observations. If $\mu_t$ is
differentiable, using the Taylor-series expansion it can be
locally approximated by a polynomial of degree $p$ of the time
distance, $j$, between $y_t$ and the neighbouring observations
$y_{t+j}$. Hence, $\mu_{t+j} \approx m_{t+j}$, with
$$m_{t+j} = \beta_0+\beta_1  j+\cdots+\beta_p j^p, j = 0, \pm 1, \ldots,\pm h.$$
The degree of the polynomial is crucial in determining the
accuracy of the approximation. Another essential quantity is the
size $h$  of the neighbourhood around time $t$; for
$t=h+1,...,n-h+1$, the neighbourhood consists of $2h+1$
consecutive and regularly spaced time points at which observations
$y_{t+j}$ are made. At the boundaries, asymmetric neighborhood
will be considered. The parameter $h$ is the bandwidth, for which
we assume $p \leq 2h$ throughout.

Replacing  $\mu_{t+j}$ by its approximation gives the local
polynomial model:
\begin{equation}
y_{t+j} =  \sum_{k=0 }^p \beta_k j^k+\epsilon_{t+j}, \;\;\;j = 0,
\pm 1, \ldots,\pm h. \label{eq:lam}
\end{equation}
Assuming that $\epsilon_{t+j}\sim\NID(0,\sigma^2)$, then
(\ref{eq:lam}) is a linear Gaussian regression model with
explanatory variables given by the powers of the time distance
$j^k, k= 0, \ldots, p$ and unknown coefficients $\beta_k$, which
are proportional to the $k$-th order derivatives of $\mu_t$.
Working with the linear Gaussian approximating model, we are faced
with the problem of estimating $m_{t}= \beta_0$, i.e. the value of
the approximating polynomial for $j=0$, which is the intercept
$\beta_0$ of the approximating polynomial.

The model  (\ref{eq:lam}) can be rewritten in matrix  notation as
follows: \begin{equation}\by = \bX \bbeta + \bepsilon, \;\;\;
\bepsilon \sim \N(\0,\sigma^2\bI)\label{eq:model}\end{equation}
where $\by= \left[y_{t-h}, \cdots, y_{t}, \cdots,
y_{t+h}\right]'$, $\bepsilon= \left[\epsilon_{t-h}, \cdots,
\epsilon_{t},  \cdots, \epsilon_{t+h}\right]'$,
$$\bX=\left[\begin{array}{cccc}
1 & -h &  \cdots & (-h)^p \\
1 & -(h-1)  & \cdots & [-(h-1)]^p \\
\vdots  & \vdots & \cdots & \vdots \\
1 & 0   & \cdots & 0 \\
\vdots  & \vdots & \cdots & \vdots \\
1 & h-1 & \cdots & (h-1)^p \\
1 & h   & \cdots & h^p\\
\end{array}\right],
\bbeta =  \left[\begin{array}{l} \beta_0 \\  \beta_1  \\ \vdots \\
\beta_p\end{array}\right].$$

Provided that $p\leq 2h$, the $p+1$ unknown coefficients $\beta_k,
k=0,\ldots,p,$ can be estimated by the method of weighted least
squares which consists of minimising with respect to the
$\beta_k$'s the objective function:
$$S(\hat{\beta}_0, \ldots,\hat{\beta}_p ) = \sum\limits_{j=-h}^h
\kappa_j \left( y_{t+j}-\hat{\beta}_0-\hat{\beta}_1
j-\hat{\beta}_2 j^2 -\cdots-\hat{\beta}_p j^p\right)^2,$$ where
$\kappa_j\geq 0$ is a set of weights that define, either
explicitly or implicitly, a kernel function. In general, kernels
are chosen to be symmetric and non increasing functions of $j$, in
order to weight the observations differently according  to  their
distance from time $t$; in particular, larger weight may be
assigned to the observations that are closer to $t$. As a result,
the influence of each individual observation is controlled not
only by the bandwidth $h$ but also by the kernel. Defining  $\bK =
\diag(\kappa_h, \ldots, \kappa_1, \kappa_0, \kappa_1, \ldots,
\kappa_h)$, the WLS estimate of the coefficients is
$\hat{\bbeta}=(\bX'\bK\bX)^{-1}\bX'\bK\by$. In order to obtain
$\hat{m}_t= \hat{\beta}_0$, we need to select the first element of
the vector $\hat{\bbeta}$. Hence, denoting by $\be_1$ the $p+1$
vector $\be_1' = [1, 0, \ldots, 0]$,
$$\hat{m}_t =  \be_1'\hat{\bbeta}= \be_1'(\bX'\bK\bX)^{-1}\bX'\bK\by = \bw'\by = \sum_{j=-h}^h\w_j y_{t-j},$$
which expresses the estimate of the trend as  a linear combination
of the observations with coefficients
\begin{equation}\bw' = \be_1'(\bX'\bK\bX)^{-1}\bX'\bK.
\label{eq:w}
\end{equation}

The  linear combination yielding  the trend estimate is the local
polynomial two-sided  filter. It satisfies $\bX' \bw = \be_1$. As
a consequence, the filter $\bw$ is said to   preserve a
deterministic polynomial  of order $p$. Moreover, the filter
weights are symmetric ($\w_j=\w_{-j}$), which follows from the
symmetry of the kernel weights $\kappa_j$, and the assumption that
the available observations are equally spaced.

As an example that we shall adopt in the following, we consider
the Henderson filter (Henderson, 1916; see also Kenny and Durbin,
1982, Loader, 1999, Ladiray and Quenneville, 2001) that arises as
the weighted least squares estimator of a local cubic trend at
time $t$ using the kernel $\kappa_j =
[(h+1)^2-j^2][(h+2)^2-j^2][(h+3)^2-j^2]$. These weights minimise
the variance of the third differences of the estimated trend
(maximum smoothness criterion), subject to the cubic reproducing
property.

\subsection{Asymmetric filters for the estimation at the boundaries \label{sec:bound}}

The derivation of the two-sided symmetric filter has  assumed the
availability of $2h+1$ observations centred at $t$. Obviously, for
a given finite sequence $y_t, t=1,\ldots, n$, it is not possible
to obtain the estimates of the signal for the (first and) last $h$
time points, which is inconvenient, since we are typically most
interested at the most recent estimates.

We can envisage  three fundamental approaches to the estimation of
the signal at the extremes of the sample period:
\begin{enumerate}
\item the construction of asymmetric filters that result from
fitting a local polynomial to the available observations $y_t$,
$t=n-h+1, n-h+2, \ldots, n$; \item the application of the
symmetric two sided filter $\bw$ to the series extended by $h$
forecasts $\hat{y}_{n+l|n}, l=1,\ldots,h$,
  (and backcasts $\hat{y}_{1-l|n}$);
\item the derivation of the asymmetric filter which minimises the
revision mean square  error subject to polynomial reproducing
constraints.
\end{enumerate}

The  trend estimates  for the last $h$ data points,
$\hat{m}_{n-h+1|n}, \ldots, \hat{m}_{n|n}$, use  respectively $2h,
2h-1,\ldots,  h+1$ observations. It is thus inevitable that the
last $h$ estimates of the trend will be subject to revision as new
observations become available. In the sequel we shall denote by
$q$ the number of future observations available at time $t$ (the
period which our estimate is referred to), $q=0, \ldots, h$, and
by $\hat{m}_{t|t+q}$ the estimate of the signal at time $t$ using
the information available up to time $t+q$, with $0 \leq q \leq
h$; $\hat{m}_{t|t}$ is usually known as the real time estimate
since it uses only the past and current information.

We now review the first strategy, which results from the automatic
adaptation of the local polynomial filter to the available sample;
we then interpret the results in terms of the other two
strategies. The approximate model $y_{t+j} =
m_{t+j}+\epsilon_{t+j}$ is assumed to hold for
 $j=-h,-h+1, \ldots, q$, and the estimators of the coefficients
 $\hat{\beta}_k$, $k=0,\ldots, d$,  minimise
$$\begin{array}{lll}S(\hat{\beta}_0, \ldots,\hat{\beta}_d)&=&
\sum\limits_{j=-h}^{q}\kappa_j  \left(y_{t+j}-
\hat{\beta}_0-\hat{\beta}_1 j-\hat{\beta}_2
j^2-\cdots-\hat{\beta}_d j^d\right)^2.
\end{array}$$
Let us partition the matrices $\bX$, $\bK$ and the vector $\by$ as
follows:
$$\bX = \left[\begin{array}{c} \bX_p\\ \bX_f \end{array}\right],\;\;\; \by =
\left[\begin{array}{c} \by_p\\ \by_f\end{array} \right],\;\;\; \bK
= \left[\begin{array}{cc} \bK_p & \0 \\ \0 &
\bK_f\end{array}\right],$$ where $\by_p$ denotes the set of
available observations, whereas $\by_f$ is missing and $\bX$ and
$\bK$ are partitioned accordingly. The local polynomial regression
(LPR) filters arising as the solution to the above weighted least
squares problem are written in matrix notation as:
\begin{equation}
\bw_a =\bK_p\bX_p(\bX_p'\bK_p\bX_p)^{-1}\be_1. \label{eq:daf}
\end{equation} Equivalently
\begin{equation}
\bw_a=\bw_p +\bK_p\bX_p  (\bX_p'\bK_p\bX_p)^{-1} \bX_f'\bw_f,
\label{eq:asym}
\end{equation} that is obtained by partitioning the two-sided symmetric filter in two
groups, $\bw = [\bw_p', \bw_f']'$, where $\bw_p$ contains the
weights attributed to the past and current observations and
$\bw_f$ those attached to the future unavailable observations. The
proof of (\ref{eq:asym}) can be found in Proietti and Luati
(2007), where detailed proofs of other results that will be used
in this section, such as (\ref{eq:asyf}), are also available.
Equation (\ref{eq:asym}) represents the fundamental relationship
which states how the asymmetric LPR filter weights are obtained
from the symmetric ones. Premultiplying both sides by $\bX_p'$, we
can see that the asymmetric filter weights satisfy the polynomial
reproduction constraints $\bX_p'\bw_a = \bX_p'\bw_p+\bX_f'\bw_f =
\bX'\bw.$ Thus, the bias in estimating an unknown function of time
has the same order of magnitude as in the interior of time
support.

The filter resulting from the automatic adaptation of the local
polynomial fit can be equivalently derived using the second
strategy, assuming that the future observations are generated
according to a polynomial function of time of degree $p$, so that
the optimal forecasts are generated by the same polynomial model.

The third strategy consists in determining the asymmetric filter
$\bv$ minimising the mean square revision error subject to
constraints. Let us rewrite the regression model (\ref{eq:model})
as $$\by =\bU\bgamma + \bZ\bdelta + \bepsilon, \bepsilon\sim
\N(\0, \bD),$$ where we have partitioned the columns of the design
matrix $\bX=[\bU | \bZ]$, in order to separate the polynomial
constraints imposed to the filter from those assumed for the
trend. Specifically, the constraints are specified as follows:
$\bU_p' \bv =\bU' \bw,$ where $\bU = [\bU_p',\bU_f']'$. Writing
$\bD = \diag(\bD_p, \bD_f)$, the set of asymmetric weights
minimises with respect to $\bv$ the following objective function
\begin{equation}
\varphi(\bv) =
(\bv-\bw_p)'\bD_p(\bv-\bw_p)+\bw_f'\bD_f\bw_f+\left[\bdelta'(\bZ_p'\bv-\bZ'\bw)\right]^2+2\bl'
(\bU_p' \bv -\bU' \bw). \label{eq:objf}
\end{equation}
The revision error arising in estimating the signal $m_t$ is
$\hat{m}_{t|t}-\hat{m}_{t}=\bv'\by_p-\bw'\by$. Replacing $\by_p
=\bU_p\bgamma+ \bZ_p\bdelta+\bepsilon_p$, and $\by=\bU\bgamma
+\bZ\bdelta+\bepsilon$, and using $\bU_p' \bv =\bU' \bw=\0$,
 we obtain
$\hat{m}_{t|t}-\hat{m}_{t}=(\bv'\bZ_p-\bw'\bZ)\bdelta +
\bv'\bepsilon_p-\bw'\bepsilon$, where $\bepsilon=
[\bepsilon_p',\bepsilon_f']'$. Hence, the first three summands of
(\ref{eq:objf}) represent the mean square revision error, which is
broken down into the revision error variance (the first two terms)
and the squared bias term
$\left[\bdelta'(\bZ_p'\bv-\bZ'\bw)\right]^2$. The vector $\bl$ is
a vector of Lagrange multipliers. The solution is \begin{equation}
\bv = \bw_p + \bL\bU_f'\bw_f +\bM\bZ_p\bdelta\bdelta'\bZ_f'\bw_f,
\label{eq:asyf}
\end{equation}
with
$$ \bM  = \bQ^{-1}-\bQ^{-1}\bU_p[\bU_p'\bQ^{-1} \bU_p]^{-1}\bU_p'\bQ^{-1},
\bL = \bQ^{-1}\bU_p[\bU_p'\bQ^{-1} \bU_p]^{-1}.$$ The matrices
$\bM$ and $\bL$ have the following properties: $\bU_p'\bM=\0,
\bU_p'\bL = \bI.$ It should be noticed that the LPR filters arise
in the case $\bD =\bK^{-1}$ and $\bU=\bX$, so that the bias term
is zero.

The merits of the class of filters (\ref{eq:asyf}), relative to
the LPR asymmetric filters, lie in the bias-variance trade-off. In
particular, the bias can be sacrificed  for improving the variance
properties of the corresponding asymmetric filter.

\section{Matrices associated with local polynomial regression estimators} \label{sec:sm}
Any linear operator acting on an $n$-dimensional time series
$\mathbf{y}$ to produce smooth estimates of the underlying trend
can be represented in matrix form as
\begin{equation*}
\mathbf{Sy} = \hat{\bm}
\end{equation*}
where $\mathbf{S}$ is the $n\times n$ smoothing matrix
representative of a weighted average to be applied to the
observations in moving manner and $\by$ is, from now on, the
$n-$dimensional vector containing all the observations. In
practice, $\mathbf{S}$ can be constructed as the matrix
canonically associated with the linear transformation $s$, so that
its columns contain the coordinates of the $s$-transformed vectors
of the canonical basis $\mathcal E = \{\mathbf{e}_{1},
\mathbf{e}_{1}, \hdots \mathbf{e}_{n}\}$, where $\mathbf{e}_{j}$
is the vector with all zeros except for the $j$-th element, equal
to one, taken with respect to the canonical basis itself, i.e.
$\bS = \left[s(\be_1)_{\mathcal E} \cdots s(\be_n)_{\mathcal
E}\right]$.

The rows of $\mathbf{S}$, denoted by $\mathbf{w}_{t}'$, are the
filters and the generic element $\w_{tj}$ is the weight to be
assigned to the observation $y_{j}$ to get the estimated value
$\hat{m_{t}}$. The weights are null outside a bandwidth whose
length, a function of $h$, depends on the local estimation method.
In general, the $n-2h$ central values are estimated by applying
$2h+1$ symmetric weights to consecutive observations centred in
$t$ whereas the first and last $h$ trend estimates are obtained by
applying asymmetric filters of variable length to the available
observations at the boundaries of the series. Thus it follows that
$\bS$ is a banded matrix with the following structure
\begin{equation}
\mathbf{S}=\left[
\begin{array}{ccc}
\mathbf{S}_{\left( h\times 2h\right) }^{a} &  & \mathbf{O}_{\left(
h\times
n-2h\right) } \\
& \mathbf{S}_{\left( n-2h\times n\right) }^{s} &  \\
\mathbf{O}_{\left( h\times n-2h\right) } &  & \mathbf{S}_{\left(
h\times 2h\right) }^{a\ast }
\end{array}
\right] \label{eq:S}
\end{equation}
where $\mathbf{S}^{s}$ is the submatrix whose rows are the
symmetric filters, while $\mathbf{S}^{a}$ and $ \mathbf{S}^{a\ast
}$ contain the asymmetric filters to be applied to the first and
last observations, respectively; the number into parentheses
indicate the dimension of the submatrices.

$\mathbf{S}$ is centrosymmetric, in that
$\w_{tj}=\w_{n+1-t,n+1-j}$; $\mathbf{S}^{s}$ is rectangular
centrosymmetric, whereas $\mathbf{S}^{a}$\ and $\mathbf{S} ^{a\ast
}$, are one t-transform of one another, where t is a linear
transformation that consists in the
pre- and post-multiplication of a matrix by the exchange matrix $\mathbf{E}%
_{k}\in \mathbb R ^{k\times k}$ having ones on the cross diagonal
(bottom left to top right) and zeros elsewhere (Dagum and Luati,
2004). For example, t$(\bS^{a\ast})=\bE_h\bS^a\bE_{2h}$.
Centrosymmetric matrices are invariant with respect to t and
preserve their structure under matrix multiplication, thus
allowing the convolution of linear filters to be a linear filter
as well. On the other hand, they are in general not symmetric,
with the consequence that their eigenvalues and eigenvectors are
complex. In dealing with real data, such as time series, this is
inconvenient. Moreover, very little is known about the analytical
form of such quantities, except that eigenvectors are either
symmetric or skew symmetric (Weaver, 1985), i.e. invariant or
equal to their opposite if premultiplied by $\bE_n$. For symmetric
matrices, some results can be found in Cantoni and Butler (1976)
and Makhoul (1981).

The rest of the paper deals with the spectral analysis of matrices
like $\bS$. In the next section, we will define the problem and
review some asymptotic results that hold in the ideal case of
doubly infinite samples. Then, the main results on the eigenvalues
and eigenvectors in finite dimension will be derived.

\section{Spectral analysis} \label{sec:sa}

The scalar $\lambda$ is an eigenvalue of $\mathbf{S}$ if there
exists a non null vector $\mathbf{x}$ such that
$\mathbf{Sx}=\lambda\bx$ and $\bx$ is the eigenvector of
$\mathbf{S}$ corresponding to $\lambda$. If we could virtually
take an infinite time series and apply the two-sided symmetric
filter to all the observations, then we would have an infinite
smoothing matrix structured like a symmetric banded Toeplitz
(SBT), with real eigenvalues and eigenvectors. Let us suppose that
the eigenvalues can be ordered in a numerable decreasing sequence,
$\lambda_1\geq \lambda_2\geq ...\geq\lambda_n\geq ...$. Hence, the
eigenvectors $\mathbf {x}_1, \mathbf {x}_2,...,\mathbf {x}_n,...,
$ can be interpreted as time series that the filter expands,
$\lambda_i> 1$, leaves unchanged, $\lambda_i= 1$, shrinks,
$\lambda_i < 1$, or suppresses, $\lambda_i=0, i=1,2,...,n,...$ .
We may ask how do these series behave and how are they modified by
the corresponding eigenvalues.

Because of their symmetric or skew symmetric nature, the
eigenvectors are likely to be interpreted as polynomials or as
periodic components. Thus, since we are dealing with matrices
associated with trend filters, what we expect is that low
frequency components associated with smooth variations of the
underlying process are represented by long period eigenvectors and
associated with eigenvalues close to unity. On the other hand, we
expect that high frequency components associated with erratic
fluctuations will be represented by short period eigenvectors
associated with eigenvalues close to zero. Hence the eigenvectors
of $\bS$ can be interpreted as the periodic latent components of
any time series, modified by the filter through multiplication by
the corresponding eigenvalues. In fact, let us consider the linear
combination
$$\mathbf {y} = \alpha_1\mathbf {x}_1 + \alpha_2 \mathbf {x}_2
+...+ \alpha_n\mathbf {x}_n + ... $$ then
$$\mathbf{Sy}=\sum_{i=1}^{k} \lambda_i\alpha_i\mathbf {x}_i
+\sum_{i=k+1}^{\infty} \lambda_i\alpha_i \mathbf {x}_i,$$ where
the $\alpha$ depend on the series $\by$, in that they re-scale the
amplitude of each periodic component, and the $\lambda$ depend on
the smoothing matrix $\bS$, i.e. on the filter. It follows that,
independently of the $\alpha$, there will be $k$ components that
the filter leaves unchanged or smoothly shrinks, and these account
for the signal, and $\infty-k$ components that will be almost
suppressed, and these account for the noise.

The choice of $k$ turns out to be a filter design problem in time
domain. There is a mathematically elegant exact solution, which
occurs if $\mbox{rank}(\mathbf{S})=k$ that is $\hat\bm$ belongs to
the column space $\mathcal{C}(\bS)$ and $\bepsilon$ lies in the
null space $\mathcal{N}(\bS)$. In practice, even if many of the
eigenvalues are close to zero, $\bS$ is full rank and therefore we
may only look for an approximate solution that consists of
choosing a cut-off time or a cut-off eigenvalue. To do this, it is
necessary to know the analytical form, at least with some
approximations or restrictions, of the eigenvalues and
eigenvectors of $\bS$.

\subsection{Infinite dimension} \label{sec:as}

In the ideal case of a doubly infinite sample, the matrix $\bS$ is
a SBT operator whose non null elements are the Fourier
coefficients of the trigonometric polynomial (the symbol of the
matrix, see Grenander and Szeg\"{o}, 1958) $$
H(\nu)=\sum_{d=-h}^h\w_d e^{\imath\nu d}$$ and
$$\lim_{n\rightarrow\infty}\frac{1}{n}\sum_{i=1}^n
\lambda_i=\frac{1}{2\pi}\int_0^{2\pi} H(\nu)d\nu$$ with
$$\lambda_1\leq \max H(\nu), \lambda_n\geq \min H(\nu).$$
$H(\nu)$ is the transfer function of the filter evaluated at the
frequency $\nu$, expressed in radians. The fundamental eigenvalue
distribution theorem states that the spectrum of an infinite SBT
matrix is dense on the set of values that the transfer function of
the symmetric filter can assume and no revisions or phase shifts
intervene in the estimation process.

In finite dimension, the analytical form of eigenvalues and
eigenvectors is known only for few classes of matrices, which are
the tridiagonal SBT and matrices belonging to some algebras,
namely the Circulant, the Hartley and the generalised Tau. All
these matrix algebras are associated with discrete transforms such
as, respectively, the Fourier, the Hartley and the various
versions of the Sine or Cosine; see, respectively, Davis (1979),
Bini and Favati (1993), Bozzo and Di Fiore (1995) and the survey
paper by Kailath and Sayed (1995). In our setting, any algebra
undertakes different hypotheses on the future behaviour of the
series. Interpreting a smoothing matrix as the sum of a matrix
belonging to one of these algebras plus a perturbation occurring
at the boundaries, approximate results on the eigenvalues of $\bS$
can be derived. The size of the perturbation depends on the matrix
algebra and on the boundary conditions. In the following, we
consider the circulant and the reflecting algebras as well as
asymmetric filters that approximate a given two sided symmetric
filter according to a minimum mean square revision error criterion
subject to constraints.

\subsection{Circular boundary conditions} \label{sec:cbc}

The circularity assumption, that is the future behaviour of the
process is equal to its initial path, represents the ideal
situation when the transfer function of any asymmetric filter is
equal to that of the symmetric filter and no phase shifts affect
the process, like in the infinite case. However, the circularity
assumption has the limitation of being restrictive in the presence
of nonstationary trends.

In the sequel, given a two sided symmetric filter
$\{\w_{-h},...,\w_0,...,\w_h \}$, we will denote by $\bS$ the
$n\times n$ associated smoothing matrix, with boundary conditions
determined by approximate asymmetric filters, and by $\bW$ the
corresponding circulant matrix (Davis, 1979) structured like a
finite SBT plus circular corrections in the top-right and
bottom-left corners,
$$\bW=\sum_{d=0}^h \w_d\bC^d+\sum_{d=-h}^{-1} \w_d\bC^{n+d}$$
where $\bC$ is the circulant matrix (basis) whose first row is the
$n$-dimensional vector $[0,1,0,0,...,0]$. Note that $\bW$ is
symmetric, as follows by the symmetry of the filter weights. For a
square matrix $\bA$ we will denote its spectrum by $\sigma(\bA)$
and its $2$-norm by $\|\bA\|_2=\sqrt{\rho(\bA'\bA)}$ where
$\rho(\bA)$ is the spectral radius of $\bA$, which is the maximum
modulus of its eigenvalues. With this preliminary notation, we are
able to state the following result on the eigenvalues of a trend
filter matrix. The proof is in the appendix. \vspace{0.2cm}

\textbf{Theorem 1} \emph{Let $\bS$ be an $n\times n$ smoothing
matrix associated with the symmetric filter} $\{\w_{-h}$, $...$,
$\w_0,...,\w_h \}$, $n>2h$, \emph{and let $\bW$ be the
corresponding circulant matrix. Hence,
$\forall\lambda\in\sigma(\bS)$, $\exists i\in\{1,2,..,n\}$ such
that}
$$\left|\lambda-\left(\sum_{d=0}^h\w_d\omega^{(i-1)d}+
\sum_{d=-h}^{-1}\w_d\omega^{(i-1)(n+d)}\right)\right|\leq\delta_W$$
\emph{where $\omega=e^{-\imath\frac{2\pi}{n}}$ and $\delta_W=\|
\bS-\bW \|_2$.}\vspace{0.2cm}

\begin{figure}
\caption{Left. Transfer function of the symmetric Henderson
filter, $h=6$, $\nu\in[0,\pi]$ (line) and eigenvalues of the
associated circulant matrix $\bW$ (dots), $n=51$. Right.
Eigenvalue distributions of $\bW$ (dots) with asymmetric
Musgrave-LC (squares), QL (circles), CQ (stars), LPR filters
(pluses) filters, absolute values. \label{fig:FigCirc}}
\begin{center}
\vspace{-0.40cm} \hspace*{-1.00cm}
\scalebox{0.60}[0.30]{\rotatebox{-0}{\includegraphics{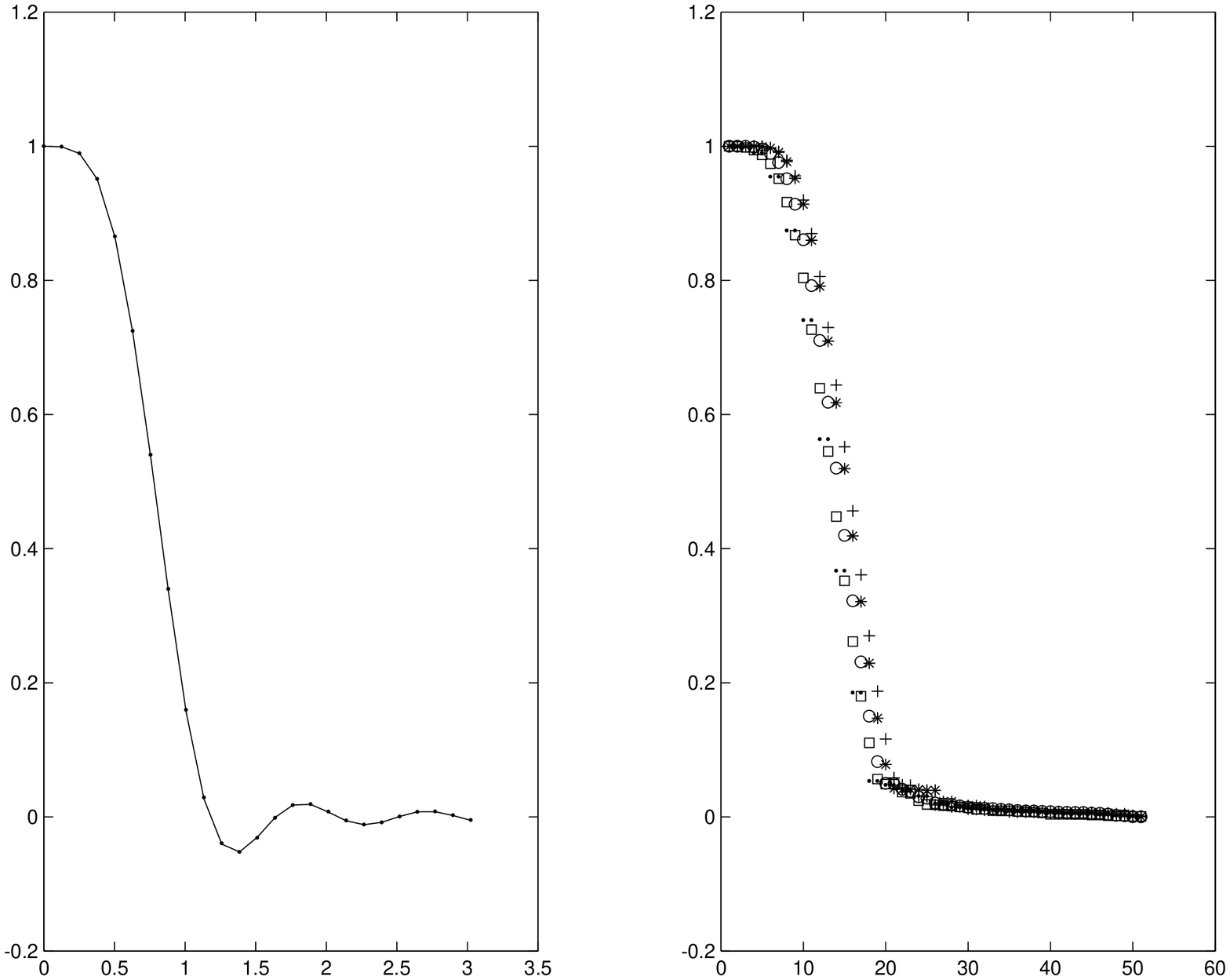}}}
\end{center}
\end{figure}

The theorem provides an upper bound on the size of the
perturbation of the eigenvalues of $\bS$ with respect to those of
$\bW$, for which an exact analytical expression is available. The
quantity $\delta_W$ measures how much the eigenvalue distribution
of a smoothing matrix moves away from that of the corresponding
circulant. On their turn, the eigenvalues of the circulant matrix
result to be distributed over the transfer function of the
symmetric filter, as the left panel of figure \ref{fig:FigCirc}
shows. What follows is that $\delta_W$ can be chosen as a measure
of how much the eigenvalue distribution of $\bS$ deviates from the
transfer function of the associated filter. In the next section,
we will show that the discrete approximation of $H(\nu)$ through
the points in $\sigma(\bW)$ can be improved by assuming the
hypothesis of reflecting behaviour of the process at the end of
the sample. As we will see, this occurs because $n$-dimensional
filtering matrices subject to reflecting boundary conditions have
$n$ distinct eigenvalues, whereas circulant matrices have pairwise
coincident eigenvalues, i.e. $\frac{n}{2}$ or $\frac{n-1}{2}+1$
distinct eigenvalues, for $n$ even or odd, respectively.

To illustrate, we consider the symmetric 13-term Henderson filter
introduced in section \ref{sec:lpr} and, as an approximation at
the boundaries, the LPR estimators and the following asymmetric
filters based on a minimum mean square revision error strategy,
subject to polynomial constraints:
\begin{description}
\item[] Linear trend - Constant fit (LC): the asymmetric LC
filters arise as the best approximations to the two-sided
Henderson filter assuming that $y_t$ is linear and imposing the
constraint that the weights sum to 1. Hence $\bU=\bi$, the unit
vector. This class contains the well-known Musgrave (1964)
surrogate filters that are commonly used to approximate the
Henderson filters. \item[] Quadratic trend - Linear fit (QL): the
asymmetric QL filters arise as the best approximations to the
two-sided Henderson filter assuming that $y_t$ is quadratic and
imposing the constraint that the estimates are capable of
reproducing a first degree polynomial. Hence $\bU$ is made of the
first two columns of $\bX$ whereas $\bZ$ contains the third column
of $\bX$. \item[] Cubic trend - Quadratic fit (CQ): the asymmetric
CQ filters arise as the best approximations to the two-sided
Henderson filter assuming that $y_t$ is a cubic function of time
and imposing the constraint that the estimates are capable of
reproducing a second degree polynomial. Hence $\bU$ is made of the
first three columns of $\bX$ whereas $\bZ$ contains the fourth
column of $\bX$.
\end{description}
Except for the LPR filters, all of the asymmetric filters are
derived here for fixed values of the parameters they depend upon,
i.e. $\delta_r^2/\sigma^2$, $r=1,2,3$ for LC, QL and CQ
respectively, and are posed equal to the value that gives the
Musgrave filter approximating the 13-term Henderson filter, i.e.
$\delta_1^2/\sigma^2=4/(3.5^2\pi)$. The parameters $\delta_1,
\delta_2, \delta_3$ represent the slope, curvature, and relative
inflexion of the trend.

The results are the following: the size of the perturbation is
minimum for the Musgrave-LC filters, being $\delta_W=0.5835$, and
maximum for the LPR filters, for which $\delta_W=1.0047$. In the
middle, the asymmetric QL, $\delta_W=0.8641$ and CQ,
$\delta_W=0.9876$. As a consequence, the eigenvalue distributions
turn out to be slight translations (towards the right) of the
absolute transfer function (the gain) of the symmetric filter:
this implies an increase in the overall variance of the estimated
trend, the increase being greater as long as $\delta_W$ increases,
as the right panel of figure \ref{fig:FigCirc} shows.

The size of the perturbation does not depend on $n$, in that the
$n-2h$ central rows of the matrix $\bS-\bW$ are all null. On the
other hand, it is highly influenced by the real time filter (last
row of $\bS$), applied to estimate the trend at time $t$ using the
available observations up to and including $t$. The fact is that,
in general, there is a discontinuity in the behaviour of the real
time filter with respect to the preceding asymmetric ones, due to
the rapid increase of the leverage of the filter, i.e. the weight
attached to the observation taken at the same time we are
estimating the trend, as long as the span of the filter decreases.
The leverage further tends to increase (up to unity) with high
degrees of the fitting polynomial (for a formal proof, see
Proietti and Luati, 2007). Here, we verify this phenomenon by
choosing as smoothing matrix the circulant matrix with first and
last rows replaced by any real time filter of the class introduced
above. The resulting values of $\delta_W$ are almost identical to
those obtained when the smoothing matrices with the whole
asymmetric filters were considered: for Musgrave-LC it is
$\delta_W=0.5247$, for QL it is $\delta_W=0.8024$, for CQ it is
$\delta_W=0.9393$, for the LPR filters  it is $\delta_W=0.9547$.
Conversely, all the values of $\delta_W$ result greater than
$0.95$ provided that the first and last row of $\bS$ are replaced
by the real time LPR filters, whose leverage is close to one.

Another factor that highly affects the size of the perturbation
(and the overall variance of the trend estimates) is the algebraic
multiplicity of the eigenvalue $\lambda=1$, that we now show to
depend on the degree of the polynomial that the filter is capable
of reproducing. The $p-$th degree polynomial reproduction
constraints met in section \ref{sec:lpr} can be written as
\begin{equation}
\langle\bw_t,\bi\rangle= 1, \;\;\;\;\;\;
\langle\bw_t,\mathbf{d}_{q}^{r}\rangle=0 \;\;\;\;\;\; \forall
t=1,2,...,n \label{eq:rp}
\end{equation} where $\bw_t'$ is
the $t$-th row of $\bS$, $\bi$ is an $n$-dimensional vector of
ones and $\mathbf{d}_{q}^{r}=[(-q)^r,(-q+1)^r,...,(n-q-1)^r]',$
with $q=t-1, \mbox{ for } t=1,...,n,   r=1,...,p$. As an example,
consider a polynomial trend $\mu_t=a_0+a_1+a_2^2+...+a_p^p$ and a
symmetric filter $\{\w_{-h},...,\w_0,...,\w_h \}$. Then
$\hat{\mu}_{t}=\sum_{d=-h}^h \w_d
[a_0+a_1(t+d)+a_2(t+d)^2+...+a_p(t+d)^p]=\mu_t $ if $\sum_{d=-h}^h
\w_d=1$ and $\sum_{d=-h}^h d^r \w_d=0$ for $r= 1,2,\ldots, p.$ The
conditions (\ref{eq:rp}) imply that
$$\bS\bi=\bi , \;\;\;\;\;\; \bS\bx_r = \bx_r$$ where $\bx_r$ is
the vector whose $t$-th coordinate is $a_0+a_1+a_2t^2+...+a_pt^r$,
that means that $\bi$ and $\bx_r, r=1,2,...,p$, are eigenvectors
of $\bS$ corresponding to the eigenvalue $\lambda =1$ of algebraic
multiplicity equal to $p+1$. It is therefore evident that the
greater the algebraic multiplicity of the eigenvalue equal to one,
the greater the displacement between the gain function of the
filter (equal to one for $\nu=0$ only) and the absolute eigenvalue
distribution.

\subsection{Reflecting boundary conditions} \label{sec:rbc}

Besides the class of circulant matrices, another class of matrices
with known spectral properties even in finite dimension is the
$\tau_{\psi\varphi}$ algebra (Bozzo and Di Fiore, 1995), that is
associated with different versions of the Sine and Cosine
trasnforms and constitutes a generalisation of the $\tau$ family
(Bini and Capovani, 1983). An $n \times n $ matrix $\bH$ belongs
to the $\tau_{\psi\varphi}$ class if and only if
$$\bT_{\psi\varphi}\bH=\bH\bT_{\psi\varphi}$$ where
$$\bT_{\psi\varphi}=\left[
\begin{array}{rrrrr}
  \psi& 1 & 0 &\cdots& 0\\
  1& 0& 1 & \ddots& 0\\
  0& 1& \ddots & \ddots& 0\\
  \vdots & \ddots&\ddots& 0& 1\\
  0 & \hdots & 0& 1 & \varphi\\
   \end{array}\right]$$
and $\psi,\varphi=0,1,-1.$ The elements $\h_{ij}$ of the matrices
in $\tau_{\psi\varphi}$ satisfy the cross sum property
$\h_{i-1j}+\h_{i+1j}=\h_{ij-1}+\h_{ij+1}$ subject to boundary
conditions determined by $\psi$ and $\varphi$. For the original
$\tau$ algebra arising when $\psi=\varphi=0$ the boundary
conditions are $\h_{0j}=\h_{i0}=\h_{n+1j}=\h_{in+1}=0$,
$i,j=1,...,n$ and all the matrices in $\tau$ can be then derived
given their first row elements. Still based on the first row of
$\bH$ but more appropriate for our purposes, since it allows to
obtain the eigenvalues and eigenvectors of $\bH$ in
$\tau_{\psi\varphi}$ in an amenable form, is the following way to
construct $\bH$  as a linear combination of powers of
$\bT_{\psi\varphi}$ (see Bini and Capovani, 1983, Proposition
2.2). Let $\bh'=[\h_1,\h_2,...,\h_n]$ be the first row of $\bH$.
Then
\begin{equation*}
\bH=\sum_{j=1}^{n}c_{j}\bT_{\psi\varphi}^{j-1}
\end{equation*} where $\bc$ is the solution of the upper
triangular system $\bQ\bc=\bh$ and $\bQ$ is the matrix whose
$j$-th column equals the first column of
$\bT_{\psi\varphi}^{j-1}$. It follows that the eigenvalues of
$\bH$ are given by
\begin{equation}
\xi_i=\sum_{j=1}^{n}\vartheta_{i}^{j-1}c_{j} \label{eq:xi}
\end{equation} where
$\vartheta_{i}, i=1,..,n,$ are the eigenvalues of
$\bT_{\psi\varphi}.$ The eigenvectors of $\bH$ are the same of
$\bT_{\psi\varphi}$.

Let us consider the reflecting hypothesis such that the first
missing observation is replaced by the last available observation,
the second missing observation is replaced by the previous to the
last observation and so on, that for a two-sided $2h+1$-term
estimator corresponds to the real time filter
$\{\w_h,\w_{h-1}+\w_h,...,\w_1+\w_2,\w_0+\w_1 \}$, made of $h+1$
terms. With the constraint of being centrosymmetric, the
reflecting matrix $\bH$ belongs to the $\tau_{11}$ algebra and its
first row is the vector
\begin{equation}\bh'=\left[ \w_0+\w_1, \w_1+\w_2,
\w_2+\w_3,...,\w_{h-1}+\w_h, \w_h,0,...,0 \right].
\label{eq:firstrow}
\end{equation}
With these premises, we are able to construct $\bH\in\tau_{11}$
corresponding to the symmetric filter
$\{\w_{-h},...,\w_0,\w_1,...,\w_h\}$ and to derive the following
result where, for sake of notation, we use the Pochhammer symbol
$(j)_q=j(j+1)(j+2)...(j+q-1)$, for $q=0,1,...,\left\lfloor
\frac{h-j-1}{2}\right\rfloor$, the latter term denoting the
largest integer less than or equal to $\frac{h-j-1}{2}$.

\vspace{0.2cm}

\textbf{Theorem 2} \emph{Let $\bS$ be an $n\times n$ smoothing
matrix associated with the symmetric filter} $\{\w_{-h}$, $...$,
$\w_0,...,\w_h \}$, \emph{and let $\bH$ be the corresponding
matrix in $\tau_{11}$. Hence, $\forall\lambda\in\sigma(\bS)$,
$\exists i\in\{1,2,..,n\}$ such that}$$\left|\lambda
-\xi_i\right|\leq\delta_H $$ \emph{where} \begin{equation}
\xi_i=\sum_{j=1}^{h+1}\left(2\cos\frac{(i-1)\pi}{n}\right)^{j-1}
\left[\w_{j-1}+\sum_{q=0}^{\left\lfloor
\frac{h-j-1}{2}\right\rfloor}\frac{(-1)^{q+1}(j)_q}{(q+1)!}(j+2q+1)\w_{j+2q+1}
\right] \label{eq:eigH}
\end{equation} \emph{and} $\delta_H=\|
\bS-\bH\|_2$.  \vspace{0.2cm}

The proof is in section \ref{sec:app}. As by-product, theorem 2
gives the eigenvalues of $\bH\in\tau_{11}$, with first row equal
to (\ref{eq:firstrow}), as an explicit function of the filter
weights, as shown in (\ref{eq:eigH}). The corresponding
eigenvectors are known (Bozzo and Di Fiore, 1995) and given by
\begin{equation}
\bz_i=k_i \left[\cos \frac{(2j-1)(i-1)\pi}{2n} \right]_j,
j=1,2,...,n \label{eq:eigvH}
\end{equation}
with $k_i=\frac{1}{\sqrt{2}}$ for $i=1$ and $k_i=1$ for $i>1.$ The
inferential procedures that will be introduced in the following
section are based on the eigenvalues and eigenvectors given by
(\ref{eq:eigH}) and (\ref{eq:eigvH}), respectively. In the sequel,
we discuss the merit of assuming reflecting rather than circulant
boundary conditions, i.e. of basing the inference on theorem 2
rather then on theorem 1.

Indeed, there are several advantages in adopting the approximation
for $\bS$ given by $\bH\in\tau_{11}$ instead of the circulant
approximation provided by $\bW$. First, all the operators
belonging to $\tau$ algebras have real eigenvalues and
eigenvectors. All the computations related to this class can
therefore be done in real arithmetic. Secondly, the reflecting
hypothesis undertaken by the $\tau_{11}$ algebra is more
appropriate than that of a circular process when the signal is a
non stationary function of time, as is the case when we are
interested in its estimate. It should be reminded that the
estimation methods considered so far are local, so that the
boundary conditions only concern a neighborhood of the ending
observations. For fixed bandwidth methods this means that only $h,
h+1,...,2h$ observations are involved in the asymmetric filtering;
if a nearest neighbourhood approach is followed, then $2h+1$
observations will be weighted even at the extremes of the series.
Another aspect that deserves to be remarked on concerns the
absolute size of the perturbation (an overestimate of the true
distance about eigenvalues), which is smaller for reflecting than
for circulant boundary conditions, i.e. $\delta_H<\delta_W$. In
fact, in general, Circulant-to-Toeplitz corrections produce
perturbations that are not smaller than Tau-to-Toeplitz
corrections, since while $\bH$ is structured as (\ref{eq:S}), the
circulant $\bW$ has nonzero corrections in the top right and
bottom left $h\times h$ blocks. When the matrix elements are the
same, this results in a greater perturbation. Table
\ref{tab:evviva} illustrates this property for the class of
approximate filters considered before.

\begin{table}[h] \centering \caption{Values of $\delta$
for $h=6$, $\tau_{11}$ and circulant algebras, approximate
asymmetric filters. \label{tab:evviva}} \vspace{0.3cm}
\begin{tabular}{lccccccc}
&     LC &     QL &     CQ & LPR \\
$\delta_H$ (Reflecting) &   0.1608  &   0.3817    &   0.7493 &  0.8351 \\
$\delta_W$ (Circulant) &  0.5835 &   0.8641  &   0.9876 & 1.0047  \\
\end{tabular}
\end{table}

\begin{figure}
\caption{Left. Transfer function of the symmetric Henderson
filter, $h=6$, $\nu\in[0,\pi]$ (line) and eigenvalues of the
associated reflecting matrix $\bH$ (crosses), $n=51$. Right. Gain
function of the symmetric Henderson filter, $h=6$ (line) and
eigenvalue distributions of $\bS$ with asymmetric Musgrave-LC
(squares), QL (circles), CQ (stars), LPR filters  (pluses)
filters. \label{fig:FigTau}}
\begin{center}
\vspace{-0.40cm} \hspace*{-1.00cm}
\scalebox{0.60}[0.30]{\rotatebox{-0}{\includegraphics{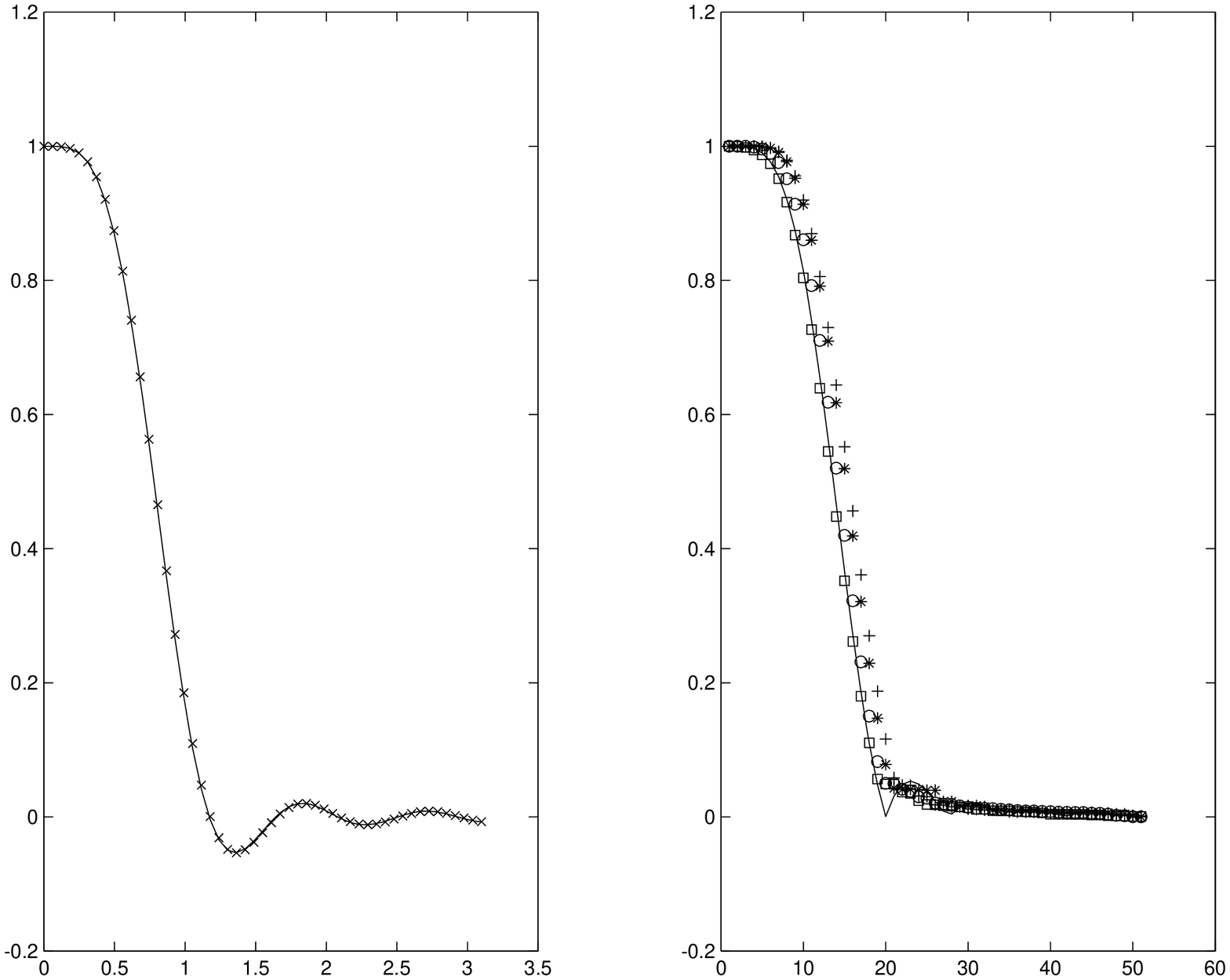}}}
\end{center}
\end{figure}

Finally, as we anticipated in the preceding subsection, the main
aspect concerning the $\tau_{11}$ approximation is that $\bH$ has
$n$ distinct eigenvalues compared to the at most $\frac{n-1}{2}+1$
of $\bW$, compare the left panel of figure \ref{fig:FigTau} with
left panel of figure \ref{fig:FigCirc}. What follows is that the
eigenvalue distribution of any smoothing matrix $\bS$ having the
form of (\ref{eq:S}) can be approximated by that of the
corresponding $\bH\in\tau_{11}$, having the same
submatrix-structure, with a deviation smaller than $\delta_H$. The
same order deviation occurs between the eigenvalue distribution of
$\bS$ and the transfer function of the corresponding symmetric
filter, as showed in the right panel of figure \ref{fig:FigTau}.

To conclude our discussion on the eigenvalues of $\bS$, we remark
that their complex part is merely generated by the finite
approximation and not related to the phase that in general affects
the asymmetric filters. This can be easily understood by means of
a counterexample: the matrix associated with a cubic smoothing
spline (see Whaba, 1990, and Green and Silverman, 1994) is
symmetric, so that its eigenvalues are real even if the asymmetric
filters do produce phase shifts.

\begin{figure}
\caption{Coordinates of the first four eigenvectors $\bz_{i}$,
$i=1,2,3,4$ of $\bH$ (crosses) and of $(\bI+\bDelta_{\bH})\bz_{i}$
(circles) plotted against $t=1,2,...n$, $n=13, h=6$, symmetric
Henderson filter. \label{fig:FigEigvector}}
\begin{center}
\vspace{-0.40cm} \hspace*{-1.00cm}
\scalebox{0.50}[0.40]{\rotatebox{-0}{\includegraphics{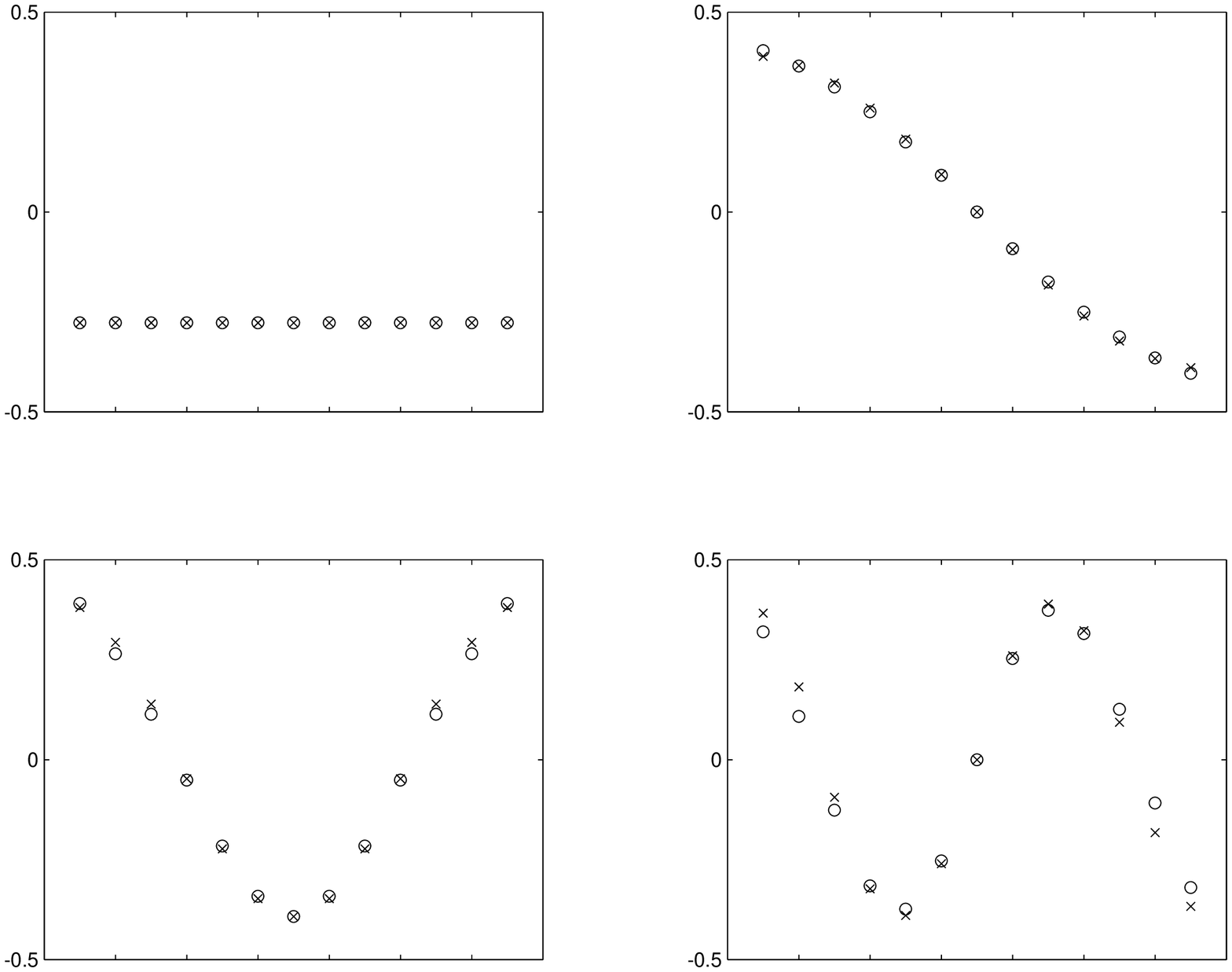}}}
\end{center}
\end{figure}

We now consider the eigenvectors. In the preceding section we have
proven that if the filter reproduces a polynomial of order $p$,
then there exist $p+1$ eigenvectors, associated with the
eigenvalue $\lambda=1$, that describe a constant ($r=0$), linear
($r=1$), quadratic ($r=2$), cubic ($r=3$) and so on up to a $p$-th
order polynomial function of the time.
In general, the analytical expression of the eigenvectors of a
smoothing matrix cannot be derived using the perturbation theory,
not even in an approximate form. However, evaluating the action of
$\bS$ on the eigenvectors of $\bH$, we are able to show that,
unless for the boundaries, the latent components of $\bS$ can be
fairly approximated by those of $\bH$. In fact, let us decompose
the time series $\by$ as a linear combination of the $n$ known
real and orthogonal latent components represented by the
eigenvectors of $\bH$,
$$\mathbf {y} = \theta_1\mathbf {z}_1 + \theta_2 \mathbf {z}_2
+...+ \theta_n\mathbf {z}_n$$ where the $\bz_i$ are given by
(\ref{eq:eigvH}) and $\btheta=[\theta_1,...,\theta_n]'$ is a
vector of coefficients. It follows from theorem 2 that
$$\bS\by=\sum_{i=1}^{n}\theta_i\xi_i\bz_i+\sum_{i=1}^{n}\theta_i\bDelta_{\bH}\bz_i$$
where $\bDelta_{\bH}\bz_i$ is a vector of zeros except for the
first and last $h$ coordinates, i.e.
$$\bDelta_{\bH}\bz_i=\left[\begin{array}{ccc}\bz_i^* \\ \mathbf{0} \\ \bE_h\bz_i^* \end{array}\right]$$
and $\bz_i^*=\sum_{j=1}^q(\bS_{ij}-\bH_{ij})\mbox{z}_{ij}$ for
$q=h+1,...,2h$ and $i=1,2,...,h$. Due to the fact that the
elements of both $\bS$ and $\bH$ add up to one and their absolute
values are in general smaller than one, the values in $\bz_i^*$
and in $\bE_h\bz_i^*$ are almost zero. This holds not only for
$n\gg h$, which is the case when we usually apply local filters,
but also for $n$ close to $h$, as figure \ref{fig:FigEigvector}
illustrates in the limiting case of $n=2h+1$, where the
approximation concerns the maximum number of boundary
approximations, namely $n-1$.

\section{Filter design in the time domain}\label{sec:fd}

The results of the preceding sections are applied for a filter
design in time domain. The aim is to obtain estimates with smaller
variance and almost equal bias than those produced by $\bS$. The
method consists of modifying $\bS$ so that $n-k$ high frequency
noisy components that the filter is not capable of eliminating are
given zero weight. This is done through the spectral decomposition
of $\bH$. The choice of $k$ i.e. of the cut-off eigenvalue $\xi_k$
will be discussed later in this section.

Decomposing $\bS=\bH+\bDelta_{\bH}$ and $\bH =\bZ\bmX\bZ'$, where
$\bmX=\diag\{\xi_1,\xi_2,...,\xi_n\}$, and writing
$\by=\bZ\btheta$, we get
$$\begin{array}{lll}
\bS\by &=& \bZ\bmX\btheta+\bDelta_{\bH}\bZ\btheta\\
 &\approx& \bZ\bmX_k\btheta+\bDelta_{\bH}\bZ\btheta
\end{array}$$
where $\bmX_k$ is the matrix obtained by replacing with zeros the
eigenvalues of $\bH$ that are smaller than a cut-off eigenvalue
$\xi_k$ and $\bDelta_{\bH}\bZ\btheta$ is a null vector except for
the first and last elements that account for the boundary
conditions. Turning to the original coordinate system and
arranging the boundaries, we get the new estimator
$$\begin{array}{lll}
\bS_k &=& \bH_k+\bDelta_k+\bDelta_{\bH}\\
 &=& \bH_{(k)}+\bDelta_{\bH}
\end{array}$$ where
$\bH_{(k)}$ is the matrix with boundaries equal to those of $\bH$
and interior equal to that of
$\bH_k=\bZ\bmX_k\bZ'$. 
In other words, $\bH_{(k)}$ is structured like (\ref{eq:S}) with
$\bH_k^a=\bH^a$, $\bH_k^{a*}=\mbox{t}(\bH^a)$  and
$\bH_k^s=[\bZ\bmX_k\bZ']^s$. Hence a new smoothing matrix is
obtained, $\bS_k$, and consequently new trend estimates, say
$\hat{\bm}_k$.

In practice, the procedure is much easier to apply. In fact, given
a symmetric filter, it consists of: obtaining $\bH$, replacing it
by $\bH_k$ and then adjusting the boundaries with suitable chosen
asymmetric filters to get $\bS_k$. Besides simplicity and variance
improvement in the interior, this procedure  allows a full choice
of the set of asymmetric weights. Indeed, in the examples we shall
illustrate at the end of this section, due to the strong impact of
the real time filter respect to all the asymmetric ones, we will
replace only the last row of $\bH_k$.

\subsection{Bias-variance trade off}\label{sec:bv}

Let us assume that $\bepsilon\sim\mathcal{N}(0,\sigma^2\bI)$. The
variance of the estimates obtained by $\bS$ is given by
$V(\hat\bm)=V(\bS\by)=\bS\sigma^2\bI\bS'
$.   
It follows that 
$$V(\hat\bm)-V(\hat\bm_k)=\sigma^2[\bZ(\bmX^2-\bmX_k^2)\bZ'+
(\bH-\bH_k)\bDelta_{\bH}'+\bDelta_{\bH}(\bH-\bH_k)'-(\bH_k\bDelta_k'+\bDelta_k\bH_k')]$$
where $\bZ(\bmX^2-\bmX_k^2)\bZ'$ is the main contribution to the
variance in the interior and is greater than zero in the sense of
a positive definite matrix; the two summands left restitute a
matrix with non null first and last $2h$ rows only, given that
$\bDelta_{\bH}$ and $\bDelta_k$ have top left and bottom right
nonzero blocks of dimension $h\times 2h$. So, even if they mainly
account for the variance at the boundaries, they also contribute
to the variance in the interior. However, for $h\ll n$ the
contribution is negligible with respect to that of the first
summand and it is common to both $\hat\bm$ and $\hat\bm_k$.

The bias is given by $B(\hat{\bm})=\bmu-E(\hat{\bm})=
\bmu-(\bH+\bDelta_{\bH})\bmu=[\bI-(\bH+\bDelta_{\bH})]\bmu.$ As
introduced by the filtering procedure, the bias is smaller as long
as $\bS$ tends to the identity matrix (in terms of the eigenvalues
of $\bS=\bI$, there are $n$ eigenvalues equal to one and therefore
the filter is capable of reproducing an $n$-degree polynomial
interpolating the data, i.e. the series itself). Comparing the
bias of the two estimators we see that
$$B(\hat\bm)-B(\hat\bm_k)=\left[\bZ(\bmX_k-\bmX)\bZ'+\bDelta_k\right]\bmu$$ and
so a measure of the discrepancy between the bias of $\hat\bm_k$
and that of $\hat\bm$, in the interior, is
$$\frac{1}{n}\mbox{ tr }\{\bmX_k - \bmX\}=-\frac{1}{n}\sum_{i=n-k+1}^n
\xi_i.$$ In general $\mbox{ tr }\{\bmX_k - \bmX\}$ is a negative
quantity that normalised by $n$ is negligible, given that the last
$n-k$ eigenvalues ar almost zero, as follows by (\ref{eq:eigH}).

The choice of $k$ is a further balancing of the trade-off between
bias and variance of the filter. The trend in the interior is made
smoother without sensibly increasing the bias. There are several
options regarding how to choose $k$. One of them, which we shall
adopt in our illustrations, consists of selecting $k$ or
equivalently $\xi_k$ that minimises the distance of the eigenvalue
distribution of $\bH$ with that of the ideal low pass filter
having first $k$ eigenvalues equal to one and and last $n-k$ equal
to zero. In other words, we look for $k$ such that
\begin{equation}f(k)=\|\bi_{(k)} - \bxi\|^2 \label{eq:min}
\end{equation} is minimum, where $\bi_{(k)}$ is an $n$ dimensional vector with first $k$
coordinates equal to one and the remaining equal to zero, whereas
$\bxi=[\xi_1, \xi_2,...,\xi_n]'$ and the $\xi_i$ are given by
(\ref{eq:eigH}). The function $f(k)=\sum_{i=1}^k(1-\xi)^2-\xi_i^2$
can be written as
$f(k)=f(k-1)+(1-\xi_k)^2-\xi_k^2=f(k-1)+(1-2\xi_k)$ and therefore
reaches its minimum for $\xi_k=0.5 .$ This strategy is equivalent
to finding the cut-off frequency that minimises the distance
between the transfer functions of the symmetric filter and of the
ideal low-pass filter $$\int_{-\pi}^\pi |I(\nu)-H(\nu)|^2d\nu $$
where $I(\nu)=1$ for $\nu<\frac{2\pi}{p}$and $I(\nu)=0$ otherwise.
The equivalence is based on the relation between time and
frequency domain. In fact, for a fixed a cut-off frequency
$\nu=\nu_k$, the cut-off time $k=\frac{\nu_k n}{\pi}$ is obtained
with a precision that increases as long as $n$ is large. For
instance, if we are given monthly data and are interested in
removing 10-month cycles that can be wrongly interpreted as
turning points of the trend curve, then we may replace by zeros
all the eigenvalues smaller than $\xi_k$ with $k=\frac{2n}{10}.$

Finally, we would like to remark that whenever the interest is in
the smoothness of the new estimator rather than in the exact value
of $k$, a graphical inspection method may be appropriate. Having
plotted the eigenvalue distribution, a suitable cut-off eigenvalue
may be directly viewed. If the choice of $k$ is not related to
formal inferential procedure (e.g. restrictions on the bias) this
method works well.

\subsection{Empirical analysis}\label{sec:ill}

In this section we provide illustrations of the eigenvalue-based
method for reducing the variance of the trend estimates obtained
with a given symmetric filter that is applied to real data. As a
symmetric estimator, we consider the $13$ term Henderson (1916)
filter, which plays a prominent role in empirical applications,
especially for trend estimation within the X-11 filter, which is
an integral part of the X-12-ARIMA procedure, the official
seasonal adjustment procedure in the U.S., Canada, the U.K. and
many other countries. See Dagum (1980), Findley {\em et al.}
(1998) and Ladiray and Quenneville (2001) for more details. As for
the asymmetric filters, the reflecting have been chosen except for
the case of the real time filter. In particular, the QL (Proietti
and Luati, 2007) and Musgrave (1964) real time filters discussed
in section \ref{sec:cbc} have been applied and compared. The
smoothing matrix $\bS$ is therefore equal to $\bH$ except for the
last row changed. To obtain $\bS_k$ we find the spectral
decomposition of $\bH$ and select the cut-off eigenvalue according
to (\ref{eq:min}), i.e. $\xi_k=0.5$.

\begin{figure}
\caption{Index of Industrial Production, Italy. Source: Istat. Top
and center. Original series with estimates obtained by $\bH$
(dotted line) and by $\bH^*$(continuous line). Bottom. Real time
estimates by QL (continuous line), reflecting (dotted line) and
Musgrave (dots) real time filters.}
\begin{center}
\vspace{-0.40cm} \hspace*{-1.00cm}
\scalebox{0.60}[0.90]{\rotatebox{-0}{\includegraphics{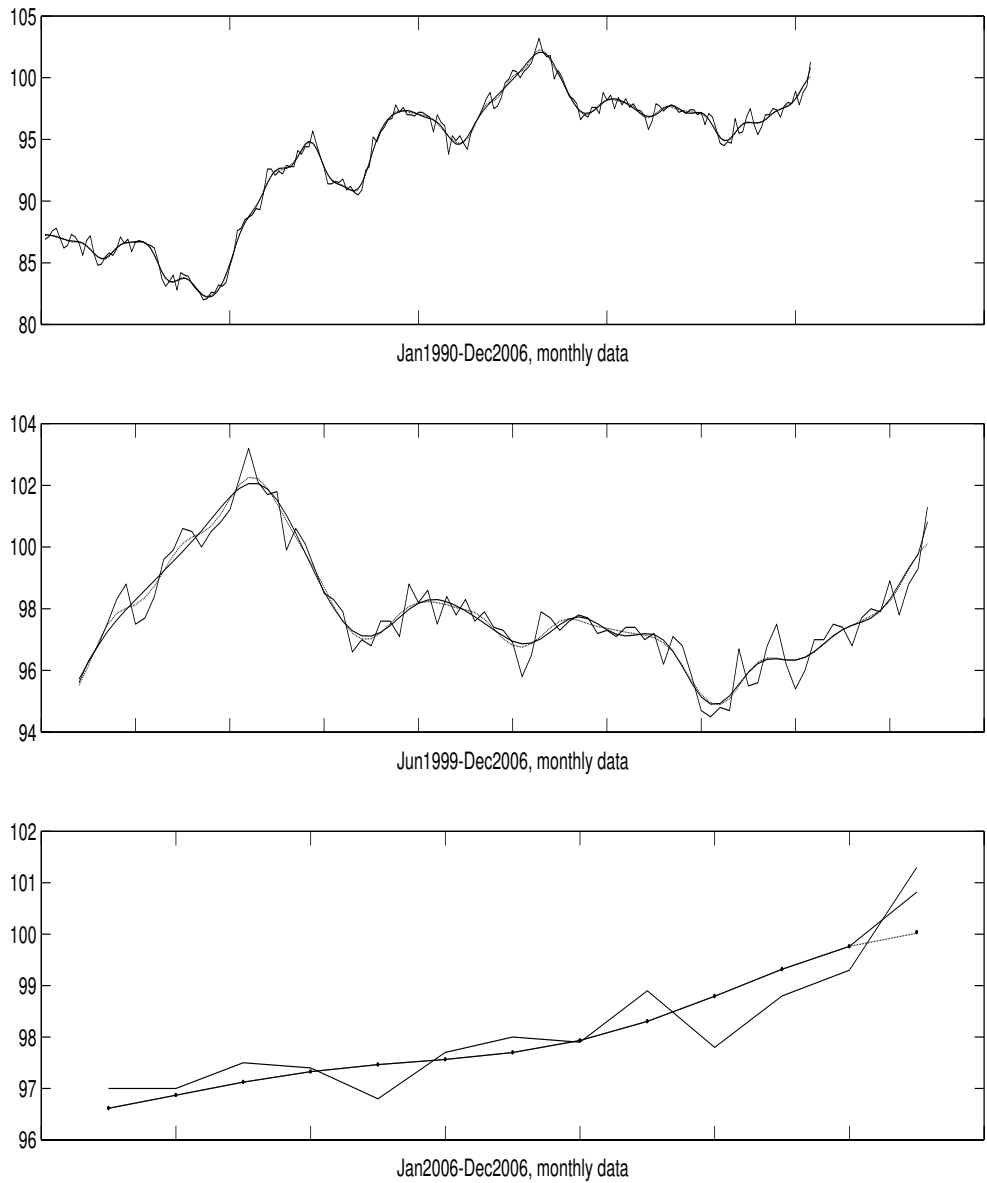}}}
\label{fig:FigIpi}
\end{center}
\end{figure}

Our first illustration deals with the Italian index of industrial
production. The top panel of figure \ref{fig:FigIpi} represents
the original series with the trend estimates $\hat\bm$ (dotted
line) and $\hat\bm_k$ (continuous line). The gain in smoothness
obtained using the latter estimator is not so evident in the whole
series but can be clearly seen in the central panel of figure
\ref{fig:FigIpi}, where a subset of the data is represented. The
estimates obtained by $\hat\bm_k$ are less sensitive to the
fluctuations of the series, note in particular the behaviour in
the period ranging from June 1999 to June 2001 when the series is
slightly increasing: the original filter estimates are sensible to
noisy fluctuations that do not affect the modified version where
highly noisy components are removed instead of just smoothed. The
bottom panel shows the last year estimates to give an idea of a
comparison among asymmetric filters: the Musgrave real time filter
(dots) behaves almost like the reflecting one (dotted line)
whereas the QL real time filter (continuous line) follows the
series increase.

\begin{figure}
\caption{Euro Area Industry, Retail Ea4. Source: European
commission. Top and center. Original series with estimates
obtained by $\bH$ (dotted line) and by $\bH^*$(continuous line).
Bottom. Real time estimates by QL (continuous line), reflecting
(dotted line) and Musgrave (dots) real time filters.}
\begin{center}
\vspace{-0.40cm} \hspace*{-1.00cm}
\scalebox{0.60}[0.90]{\rotatebox{-0}{\includegraphics{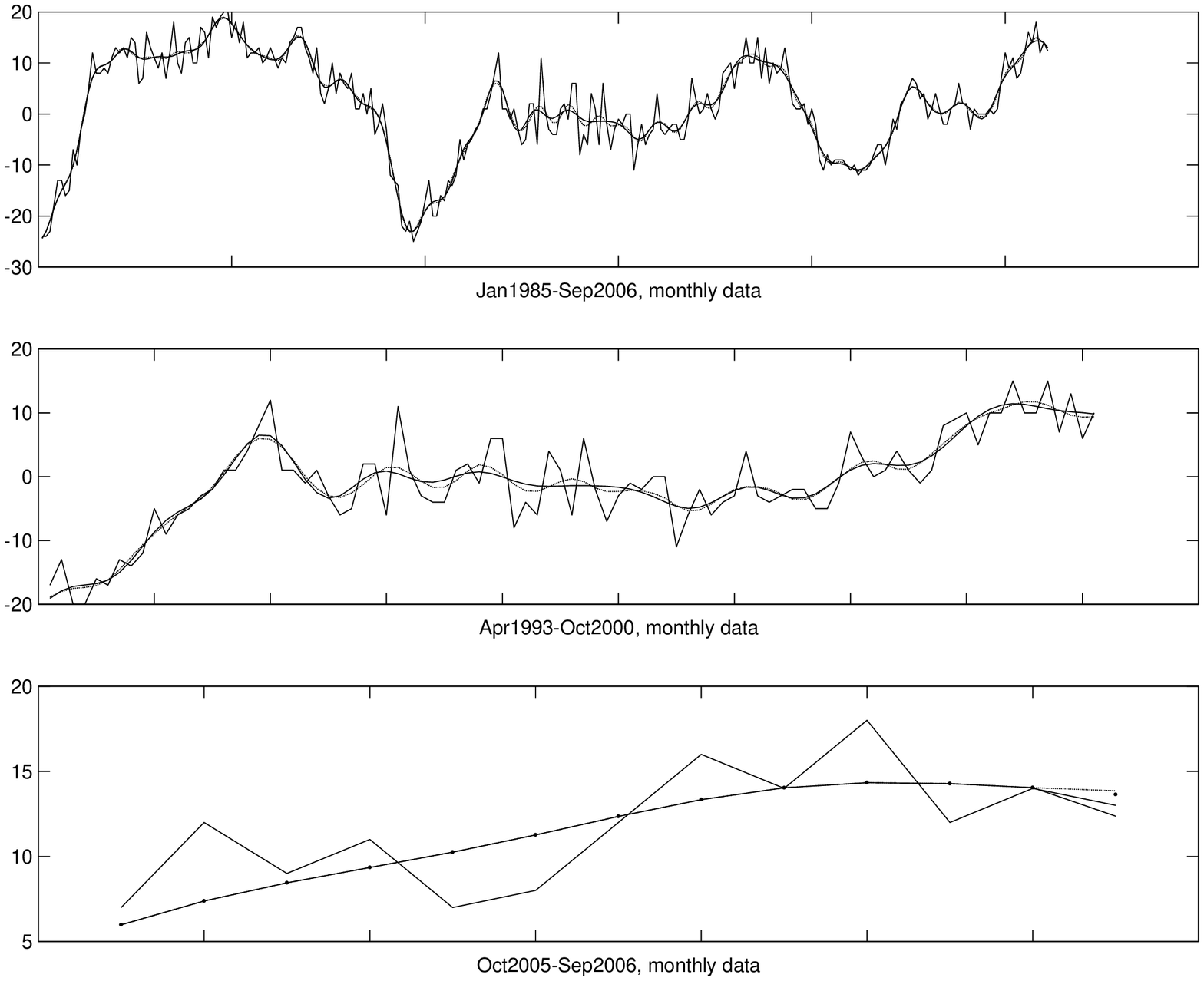}}}
\label{fig:FigRetEa4}
\end{center}
\end{figure}

Analog considerations apply to our second illustration, which
concerns the series of retails of Euro area 4, see figure
\ref{fig:FigRetEa4}. This series is affected by an increase in
variability even during periods of stationarity of the trend, as
the top panel of the figure shows. The 13-term Henderson filter
(dotted line) is known to be particularly reacting to short cycles
that, if not smoothed enough, can be falsely interpreted as false
turning points. The central part of figure \ref{fig:FigRetEa4}
illustrates that the modified estimator (continuous line) where
eigenvalues smaller than 0.5 are replaced by zeros produces
smoother trend values without affecting the capability of catching
true turning points, such as that occurred in November 1994. As in
the previous case, in the bottom panel of figure
\ref{fig:FigRetEa4} the last year estimates obtained with Musgrave
(dots), the reflecting (dotted line) and the QL (continuous line)
real time filters are illustrated. Even in this case, the QL
reacts to the changes in the direction of the series more than the
other two estimators which behave almost in the same manner.

\section{Concluding remarks} \label{sec:c}

This paper provided a decomposition of time series in periodic
latent components that depends on the underlying trend estimation
method. In particular, given a symmetric local polynomial
regression estimator with reflecting boundary conditions, the
latent components are given exactly by equation (\ref{eq:eigvH}).
These will be smoothed by an amount equal to (\ref{eq:eigH}). If
different asymmetric filters for current trend estimation are
adapted at the boundaries, then an approximation whose size was
given in theorem 2 occurs in the eigenvalue distribution.

Concerning the latter, it was shown in the paper that, in finite
dimension, an approximated version of the fundamental eigenvalue
distribution theorem holds. In fact, the eigenvalue distribution
of a trend filter matrix turned out to be a discrete approximation
of the the transfer function of the corresponding symmetric
filter. Once again, the size of the approximation depends on
boundary conditions. Circular and reflecting boundary conditions
were illustrated and discussed. In any case, it emerged that as
long as the locally weighted regression method is capable of
reproducing a high degree polynomial trend, the approximation to
the transfer function of the filter becomes worse, essentially due
to the exploding behavior of the real time filter when the degree
of the fitting polynomial increases.

It followed that, as well as the transfer function, the eigenvalue
distribution represents a measure of the overall variance left in
the input series by the smoothing procedure. More relevant, the
decomposition in periodic latent component to which smoothing
eigenvalues are associated constituted a framework for reducing
the variance of the estimates.

In fact, based on the analytical knowledge of the eigenvectors and
eigenvalues, it was possible to improve the inferential properties
of a given filter by annihilating noisy components that would have
been otherwise only smoothed. The selection of a cut-off
eigenvalue after which all the components received zero weight was
discussed and new filters with smaller variance and almost equal
bias to the original one were so derived. Applications to real
data showed the variance improvement, especially for what concerns
short cycles that may wrongly be interpreted as turning points of
the trend-cycle.

\section{Appendix} \label{sec:app}

\hspace{0.5cm}\textbf{Proof of theorem 1} The matrix $\bS$ can be
written as $\bS=\bW+\bDelta_{\bW}$, where $\bDelta_{\bW}=\bS-\bW$.
The circulant matrix $\bW$ is diagonalised by the Fourier matrix
$$\bOmega =\frac{1}{\sqrt{n}}[\omega^{(i-1)(j-1)}]_{ij},$$
$i,j=1,...,n$, satisfying $\|\bOmega\|_2\|\bOmega^{-1}\|_2=1$, and
its spectrum is $\sigma(\bW)= \{\zeta_1, \zeta_2, ...,\zeta_n\}$,
with
$$\zeta_i=\sum_{d=0}^h\w_d\omega^{(i-1)d}+\sum_{d=-h}^{-1}\w_d\omega^{(i-1)(n+d)}.$$
Setting $\delta_W=\|\bDelta_{\bW}\|_2$, the thesis follows from
the Bauer-Fike perturbation theorem applied choosing the $2$-norm
as an absolute norm (Bauer and Fike, 1960).
$\blacksquare$\vspace{0.3cm}

\textbf{Proof of theorem 2} Let us write $\bS=\bH+\bDelta_{\bH}$.
The first part of the proof is analog to the proof of theorem 1,
provided that the matrix $\bH$ is diagonalised by the orthogonal
matrix
$$\bZ=\sqrt{\frac{2}{n}}\left[k_j \cos \frac{(2i-1)(j-1)\pi}{2n}
\right]_{ij}, i,j=1,2,...,n$$ where $k_j=\frac{1}{\sqrt{2}}$ for
$j=1$ and $k_j=1$ for $j>1$ which satisfies
$\|\bZ\|_2\|\bZ^{-1}\|_2=1$. The spectrum of $\bH$ is
$\sigma(\bH)= \{\xi_1, \xi_2, ...,\xi_{n}\}$, where
$$\xi_i=\sum_{j=1}^{n}\left(2\cos\frac{(i-1)\pi}{n}\right)^{j-1}c_{j}$$
which follows by (\ref{eq:xi}) and by the fact that the
eigenvalues of $\bT_{11}$ are (Bini and Capovani, 1983)
$$\vartheta_i= 2\cos\frac{(i-1)\pi}{n}.$$ Setting $\delta_H=\|\bDelta_{\bH}\|_2$
and applying the Bauer-Fike theorem with the $2$-norm as an
absolute norm gives $$\left|\lambda
-\sum_{j=1}^{n}\left(2\cos\frac{(i-1)\pi}{n}\right)^{j-1}c_{j}\right|\leq\delta_H.
$$

We now prove that $c_{j}=0$ for $j>h+1$, so that the above
summation involves just $h+1$ terms instead of $n$. It follows by
the Cramer rule that, explicitly,
$$c_{j}=\frac{\det \bQ\left[j,\bh\right]}{\det\bQ}$$ where $\bQ\left[j,\bh\right]$
is the matrix obtained replacing the $j$-th column of $\bQ$ by the
vector $\bh$. The matrix $\bQ$ is upper triangular with ones on
the diagonal so its its determinant is equal to one and since the
generic element $h_j$ of $\bh$ is null for $j>h+1$ it follows that
$\det \bQ\left[j,\bh\right]=0$ and $c_{j}$ will be null as well.

Finally, we prove that
\begin{equation}
c_j=\w_{j-1}+\sum_{q=0}^{\left\lfloor
\frac{h-j-1}{2}\right\rfloor}\frac{(-1)^{q+1}(j)_q}{(q+1)!}(j+2q+1)\w_{j+2q+1}.
\label{eq:c}
\end{equation}
This expression can be directly verified by calculating $\det
\bQ\left[j,\bh\right]$ for all $j$. Here in the following, we
prove it by induction over $j=1,...,h+1$, with $h\in\mathbb N$.
\begin{itemize}
\item{ For $j=1$,
$c_1=\w_{0}+\sum_{q=0}^{\frac{h-2}{2}}(-1)^{q+1}2\w_{2q+2}$ which
follows by $(1)_q=q!$ and by simple algebra. The linear system
$\bQ\bc=\bh$ can be written as $\bc=\bQ^{-1}(\bh_1+\bh_2)$ with
$\bh_1=[\w_0, \w_1,...,\w_h,0,...,0]'$ and $\bh_2=[\w_1,
\w_2,...,\w_h,0,...,0]'$, both $n$-dimensional vectors. Since the
first row of $\bQ^{-1}$ is the vector $[1,-1,-1,1,1,-1,-1,...]$ we
have that
$c_1=(\w_{0}+\w_{1})-(\w_{1}+\w_{2})-(\w_{2}+\w_{3})+(\w_{2}+\w_{4})+...
+(-1)^{\left\lfloor\frac{h-2}{2}\right\rfloor+1}2\w_{2\left\lfloor
\frac{h-2}{2}\right\rfloor+2}$ and therefore (\ref{eq:c}) holds
for $j=1$.}

\item{ For $j=h$, $c_h=\w_{h-1}$ as it is immediate to see given
that the summation in $q$ was defined for non negative values of
$\frac{h-j-1}{2}$. All the more so, it implies that
$c_{h+1}=\w_{h}$. Hence we have showed that (\ref{eq:c}) holds for
$j=1$ and that, if it holds for $j=h$ then it holds for $j=h+1$.
This proves that (\ref{eq:c}) is true for all $h\in \mathbb N$.
The proof of theorem 2 is therefore complete $\blacksquare$}
\end{itemize}

\vspace{0.5cm}\emph{Acknowledgements} The authors would like to
kindly thank Dario Bini for the precious suggestions and the
useful discussions. This paper has been presented at the 56th
Session of the ISI 2007 Lisbon, August 22-29, 2007 and at the
Workshop Two Days of Linear Algebra, Bologna, 6-7 March 2008. We
are especially grateful to Daniele Bertaccini Enrico Bozzo, Fabio
Di Benedetto, Marco Donatelli and Valeria Simoncini for their very
competent comments.

\section*{References}
\begin{list}{\ }{\setlength\leftmargin{1cm}}

 \item \hskip-0.5cm
Anderson, T.W. (1971), \emph { The Statistical Analysis of Time
Series}, Wiley.

 \item \hskip-0.5cm
Bauer F., Fike C. (1960), Norms and Exclusion Theorems,
\emph{Numerische Mathematik}, 2, 137-141.

 \item \hskip-0.5cm
Bini D.,  Capovani M. (1983), Spectral and computational
properties of band symmetric Toeplitz matrices, \emph{Linear
Algebra and its Applications}, 52/53, pp. 99-126 .

 \item \hskip-0.5cm
Bini D., Favati P. (1993), On a matrix algebra related to the
discrete Hartley transform. \emph{SIMAX}, 14, 2, 500-507

 \item \hskip-0.5cm
Bloomfield, P. (2000). \emph{ Fourier Analysis of Time Series. An
Introduction}, Second edition, Wiley.

 \item \hskip-0.5cm B\"{o}ttcher A.,
Grudsky, S.M. (2005), \emph{Spectral Properties of Banded Toeplitz
Matrices}, Siam.

 \item \hskip-0.5cm
Bozzo E., Di Fiore C. (1995) On the use of certain matrix algebras
associated with discrete trigonometric transforms in matrix
displacement decomposition, \emph{Siam J. Matrix anal. Appl.}, 16,
1, 312-326.

 \item \hskip-0.5cm
 Buja A., Hastie T.J., Tibshirani R.J.(1989), Linear Smoothers
and Additive Models, \emph{The Annals of Statistics}, 17, 2,
453-555

 \item \hskip-0.5cm
Cantoni A., Butler P. (1976), Eigenvalues and eigenvectors of
symmetric centrosymmetric matrices, \emph{Linear Algebra and its
Applications}, 13, 275-288.

 \item \hskip-0.5cm
Cleveland, W.S.  and  Loader, C.L. (1996).
 Smoothing by Local Regression: Principles and Methods. In W. H\"{a}rdle and
M. G. Schimek, editors, \emph{Statistical Theory and Computational
Aspects of Smoothing}, 10–49. Springer, New York.

\item \hskip-0.5cm
Dagum, E. B. (1980).  T\emph{he X-11-ARIMA
Seasonal Adjustment Method}. Statistics Canada.

 \item \hskip-0.5cm
Dagum, E.B., Luati, A. (2004). A Linear Transformation and its
Properties with Special Applications in Time Series, \emph{Linear
Algebra and its Applications}, 338, 107-117.

 \item \hskip-0.5cm
Davis, P.J., (1979). \emph{Circulant matrices}, Wiley, New York.

\item \hskip-0.5cm Findley, D.F., Monsell, B.C., Bell, W.R., Otto,
M.C., Chen B. (1998). New Capabilities and Methods of the
X12-ARIMA Seasonal Adjustment Program, \emph {Journal of Business
and Economic Statistics}, 16, 2.

 \item \hskip-0.5cm
Green P.J. and Silverman, B.V. (1994) \emph {Nonparametric
Regression and Generalized Linear Models: a Roughness Penalty
Approach}. Chapman \& Hall, London.

 \item \hskip-0.5cm
Grenander, U., Szeg\"{o}, G. (1958). \emph{Toeplitz Forms and
Their Applications}, University of California Press.

 \item \hskip-0.5cm
Gray R.M. (2006) \emph{Toeplitz and Circulant Matrices: A review},
Foundations and Trends in Communications and Information Theory,
Vol 2, Issue 3, pp 155-239.

 \item \hskip-0.5cm
Green P.J. and Silverman, B.V. (1994) \emph {Nonparametric
Regression and Generalized Linear Models: a Roughness Penalty
Approach}. Chapman \& Hall, London.

 \item \hskip-0.5cm
Hastie T.J. and Tibshirani R.J.(1990), \emph{Generalized Additive
Models}, Chapman and Hall, London

 \item \hskip-0.5cm
Henderson, R. (1916). Note on Graduation by Adjusted Average,
\emph { Transaction of the Actuarial Society of America}, 17,
43-48.

 \item \hskip-0.5cm
Henderson R. (1924). A New Method for Graduation, \emph {
Transaction of the Actuarial Society of America}, 25, 29-40.

 \item \hskip-0.5cm
Kailath T., Sayed A. (1995) Displacement Structure: Theory and
Applications, \emph{SIAM Review}, 37, 3, 297-386.

 \item \hskip-0.5cm
Kendall M. G. (1973). \emph {Time Series}, Oxford University
Press, Oxford.

 \item \hskip-0.5cm
Kendall M., Stuart, A., and Ord, J.K. (1983). \emph {The Advanced
Theory of Statistics}, Vol 3. C. Griffin.

 \item \hskip-0.5cm
Kenny P.B., and Durbin J. (1982). Local Trend Estimation and
Seasonal Adjustment of Economic and Social Time Series,
\emph{Journal of the Royal Statistical Society A}, 145, I, 1-41.

 \item \hskip-0.5cm
Ladiray, D. and Quenneville, B. (2001). \emph{ Seasonal Adjustment
with the X-11 Method} (Lecture Notes in Statistics),
Springer-Verlag, New York.

 \item \hskip-0.5cm
Jenkins, G.M. and Watts, D.G. (1968). \emph{ Spectral Analysis and
its Applications}, Holden-Day, San Francisco.

 \item \hskip-0.5cm
Ladiray, D. and Quenneville, B. (2001). \emph{ Seasonal Adjustment
with the X-11 Method} (Lecture Notes in Statistics),
Springer-Verlag, New York.

 \item \hskip-0.5cm
Loader, C. (1999). \emph {Local Regression and Likelihood}.
Springer-Verlag, New York

\item \hskip-0.5cm
Makhoul J. (1981) On the Eigenvectors of
Symmetric Toeplitz Matrices, \emph{IEEE Transaction on Acoustic
Speech and Signal Processing}, 29, 4, 868-872.

 \item \hskip-0.5cm
Musgrave, J. (1964).  A Set of End Weights to End All End Weights.
Working paper. Census Bureau, Washington.

 \item \hskip-0.5cm
Proietti, T., Luati, A. (2007). Real Time Estimation in Local
Polynomial Regression, with an Application to Trend-Cycle
Analysis, submitted.

 \item \hskip-0.5cm
Priestley, M.B. (1981), \emph{Spectral Analysis and Time Series},
 Academic Press.

 \item \hskip-0.5cm
Weaver J.R. (1985), Centrosymmetric (Cross-Symmetric) Matrices,
their Basic Properties, Eigenvalues, Eigenvectors, \emph{Amer.
Math. Monthly}, 92, 711-717.

 \item \hskip-0.5cm
Wahba G. (1990), \emph{Spline Models for Observational Data},
SIAM\ Philadelphia.

 \item \hskip-0.5cm
Whittaker E.T. (1923), On a New Method of Graduation,
\emph{Proceedings of the Edinburgh Mathematical Association}, 78,
81-89

\end{list}

\end{document}